\documentclass[a4paper]{article}
\usepackage{latexsym, amsfonts, amsmath, amssymb}
\usepackage{tikz-cd}
\usepackage{microtype}
\usepackage[utf8]{inputenc}
\newcommand{\be}{\begin{equation}}
\newcommand{\ee}{\end{equation}}
\newcommand{\bea}{\begin{eqnarray}}
\newcommand{\eea}{\end{eqnarray}}
\newcommand{\bean}{\begin{eqnarray*}}
\newcommand{\eean}{\end{eqnarray*}}

\newcommand{\brray}{\begin{array}}
\newcommand{\erray}{\end{array}}
\newcommand{\ben}{\begin{equation}{nonumber}}
\newcommand{\een}{\end{equation}{nonumber}}

\newtheorem{dfn}{Definition}[section]
\newtheorem{thm}[dfn]{Theorem}
\newtheorem{lmma}[dfn]{Lemma}
\newtheorem{ppsn}[dfn]{Proposition}
\newtheorem{crlre}[dfn]{Corollary}
\newtheorem{xmpl}[dfn]{Example}
\newtheorem{rmrk}[dfn]{Remark}
\newcommand{\bdfn}{\begin{dfn}}
\newcommand{\bthm}{\begin{thm}}
\newcommand{\blmma}{\begin{lmma}}
\newcommand{\bppsn}{\begin{ppsn}}
\newcommand{\bcrlre}{\begin{crlre}}
\newcommand{\bxmpl}{\begin{xmpl}}
\newcommand{\brmrk}{\begin{rmrk}}
\newcommand{\edfn}{\end{dfn}}
\newcommand{\ethm}{\end{thm}}
\newcommand{\elmma}{\end{lmma}}
\newcommand{\eppsn}{\end{ppsn}}
\newcommand{\ecrlre}{\end{crlre}}
\newcommand{\exmpl}{\end{xmpl}}
\newcommand{\ermrk}{\end{rmrk}}

\newcommand{\IC}{\mathbb{C}}

\newcommand{\IR}{\mathbb{R}}

\newcommand{\IT}{\mathbb{T}}

\newcommand{\IZ}{\mathbb{Z}}

\newcommand{\oneform}{{\Omega}^1 ( \mathcal{A} )}

\newcommand{\twoform}{{\Omega}^2( \mathcal{A} )}
\newcommand{\tensora}{\otimes_{\mathcal{A}}}
\newcommand{\tensorsym}{\otimes^{{\rm sym}}_{\mathcal{A}}}
\newcommand{\tensorc}{\otimes_{\mathbb{C}}}
\newcommand{\A}{\mathcal{A}}

\newcommand{\E}{\mathcal{E}}
\newcommand{\F}{\mathcal{F}}
\newcommand{\Acenter}{\mathcal{Z}( \mathcal{A} )}
\newcommand{\Ecenter}{\mathcal{Z}( \mathcal{E} )}

\newcommand{\Psym}{P_{\rm sym}}
\newcommand{\Hom}{{\rm Hom}}

\newcommand{\cla}{{\cal A}}
\newcommand{\clb}{{\cal B}}

\newcommand{\cle}{{\cal E}}

\newcommand{\clh}{{\cal H}}

\newcommand{\clz}{{\cal Z}}

\def\a*{{\cal A}_{h,*}}
\def\B{{\cal B}(h)}
\def\B1{{\cal B}_1(h)}
\def\b{{\cal B}^{\rm s.a.}(h)}
\def\b1{{\cal B}^{\rm s.a.}_1(h)}

\newcommand{\ot}{\otimes}

\newcommand{\id}{\mbox{id}}

\def \qed {$\Box$}

\def\a*{{\cal A}_{h,*}}
\def\B{{\cal B}(h)}
\def\B1{{\cal B}_1(h)}
\def\b{{\cal B}^{\rm s.a.}(h)}
\def\b1{{\cal B}^{\rm s.a.}_1(h)}

\newcommand{\ev}{\rm ev}
\newcommand{\zero}{{(0)}}
\newcommand{\one}{{(1)}}

\newcommand{\RNum}[1]{\uppercase\expandafter{\romannumeral #1\relax}}

\begin{document}
\begin{center}
{\Large{\bf Levi-Civita connections for conformally deformed metrics on tame differential calculi}}\\
\vspace{0.2in}
{\large {*Jyotishman Bhowmick, *Debashish Goswami and **Soumalya Joardar}}\\
*Indian Statistical Institute, Kolkata, India\\
** IISER Kolkata, Mohanpur, Nadia, India.\\
Emails: jyotishmanb$@$gmail.com, goswamid$@$isical.ac.in, soumalya.j$@$gmail.com \\
\end{center}
\begin{abstract}
Given a tame differential calculus  over a noncommutative algebra $\A$ and an $\A$-bilinear pseudo-Riemannian metric $g_0,$ consider the conformal deformation $ g = k. g_0, $ $k$ being an invertible element of $\A.$  We prove that there exists a unique connection $\nabla$ on the bimodule of one-forms of the differential calculus which is torsionless and compatible with $g.$ We derive a concrete formula connecting $\nabla$ and the Levi-Civita connection for the pseudo-Riemannian metric $g_0.$ As an application, we  compute the Ricci and scalar curvature for a general conformal perturbation of the canonical metric on the noncommutative $2$-torus as well as for a natural metric on the quantum Heisenberg manifold. For the latter, the scalar curvature turns out to be a negative constant. 
\end{abstract}

\section{Introduction} 

In  recent years, the study of Levi-Civita connections in noncommutative geometry has attracted a lot of attention. Probes into the existence and uniqueness of Levi-Civita connections for noncommutative manifolds go back to the works of \cite{frolich}, \cite{chak_sinha} and \cite{heckenberger_etal}. The recent surge in activities around similar questions was pioneered by Rosenberg's work ( \cite{Rosenberg} ) on the noncommutative torus where the existence of a unique Levi-Civita connection on a certain bimodule of derivations was proved. This line of attack was further pursued by Peterka and Sheu ( \cite{sheu} ), Arnlind and his collaborators ( \cite{pseudo}, \cite{tiger} and references therein ) and more recently by Landi and Arnlind ( \cite{cylinder} ).

An alternative approach to the question of existence of Levi-Civita connections on the space of forms was taken by a number of mathematicians. In particular, bimodule Levi-Civita  connections as well as $\ast$-compatibility of Levi-Civita connections were studied by Beggs, Majid and their collaborators. We refer to \cite{majid_2}, \cite{majid_1} and references therein and the book \cite{beggsmajidbook} for a comprehensive account. Following this line, existence of a Levi-Civita connection on any quantum projective space ( for the Fubini-Study metric ) has been proved in \cite{matassa}.   For $q$-deformed connections on $S^3_q,$ we refer to the work of  Landi, Arnlind and Ilwale ( \cite{landiqhom} ) while for Levi-Civita connections on finite metric spaces and graphs, we refer to Chapter 8 of \cite{beggsmajidbook} and the paper \cite{sitarz} by Sitarz et al. Finally, yet another approach to study Levi-Civita connections on quasi-commutative algebras has been initiated in \cite{weber} and \cite{aschieri}.

Tame differential calculi and metric-compatibility of connections on such calculi have been studied previously in \cite{article1}, \cite{article3} and \cite{article4}. As well-known by experts, ( see \cite{beggsmajidbook} ), if $\nabla$ is a bimodule connection on the space of one-forms $\E$ of any differential calculus, then $\nabla$ induces a connection $ \nabla_{(\E \tensora \E)^*} $ on $  ( \E \tensora \E )^*:= {\rm Hom}_{\A} ( \E \tensora \E, \A ). $ In this article, we prove that this result is true for any connection provided the differential calculus is tame.  Moreover, wre show that $\nabla$ is compatible with a pseudo-Riemannian metric $g \in {\rm Hom}_{\A} ( \E \tensora \E, \A ) $ in the sense of \cite{article3} and \cite{article1} if and only if $  \nabla_{(\E \tensora \E)^*} g = 0.$  The goal of the rest of the article is to study Levi-Civita connections on conformal deformations of bilinear pseudo-Riemannian metrics ( see Definition \ref{7thjune20} ) on tame differential calculus. The main results of \cite{article1} and \cite{article3} were to prove the existence and uniqueness of Levi-Civita connections for any bilinear pseudo-Riemannian metric on tame differential calculi. In \cite{article2}, this result was extended to the class of strongly $\sigma$-compatible pseudo-Riemannian metric ( see Definition 3.1 of \cite{article2} ) which include bilinear pseudo-Riemannian metrics and their conformal deformations. However, in this article, we present a very simple proof of the existence and uniqueness of Levi-Civita connections for such conformally deformed pseudo-Riemannian metric following a  sufficient condition established in \cite{article1}. The main benefit to our approach is that we can derive ( Theorem \ref{19thjan4} ) a formula which shows how the Levi-Civita connection deforms under the conformal deformation of the pseudo-Riemannian metric. In case the bimodule of one-forms of the tame differential calculus under question is free, we can also derive the Christoffel symbols of the deformed Levi-Civita connection. Then we follow \cite{article3} to define and study curvature of a Levi-Civita connection on tame differential calculi. We end the article with some concrete examples of curvature computation. 

We should mention that in \cite{article3}, a more direct and elegant proof of existence of Levi-Civita connections for {\bf bilinear} pseudo-Riemannian metrics was derived by imitating the classical proof and delivering a Koszul-formula for the connection. We have been unable to generalize the proof of \cite{article3} to conformally deformed metrics which are typically not bilinear. 
 
  
Now we come to the class of examples for which our result is satisfied. We should admit that the class of tame differential calculi is restrictive. Nonetheless, we have several interesting examples for which we refer to Example 2.10 of \cite{article2}. This includes the differential calculi on the fuzzy $3$ spheres and  the quantum Heisenberg manifold coming from certain spectral triples defined in \cite{frolich} and \cite{chak_sinha} respectively. We refer to Theorem 5.4 and  Theorem 6.6 of \cite{article1}  for the proofs. A differential calculus on the fuzzy $2$-sphere was also proved to be tame in Theorem 8.5 of \cite{article3}.     Theorem 3.4 of \cite{soumalya} proves that a differential calculus constructed on the Cuntz algebra ( from a natural $C^*$-dynamical system )  is also tame. 

      Another important class of examples comes from Connes-Landi deformations ( \cite{connes_landi}, \cite{Connes-dubois} ) of the classical spectral triple of compact Riemannian manifolds equipped with a free and isometric toral action. Indeed, if $M$ is such a manifold, then the $C^*$-algebra $ C ( M ) $ can be deformed ( as in \cite{rieffel} ) to a possibly noncommutative $C^*$-algebra. The prescription of \cite{connes_landi} delivers a canonical spectral triple over  a canonical dense $\ast$-subalgebra of this deformation. Theorem 7.1 of \cite{article1} verifies that the differential calculus corresponding to this spectral triple is indeed tame. 

The plan of this article is as follows: in Section \ref{formulation}, we recall the definition and properties of tame differential calculi from \cite{article4}. In particular, we have  pseudo-Riemannian metrics on such calculi. We begin Section \ref{nablag} by recalling from \cite{article3} ( and \cite{article1} ) the notion of metric-compatibility of a connection on a tame differential calculus. Then we show that if $\nabla$ is a right connection on the bimodule of one-forms $\E$ of a tame differential calculus, then it induces a left connection $ \nabla_{( \E \tensora \E  )^*} $ on the bimodule $ ( \E \tensora \E )^*:= {\rm Hom}_{\A} ( \E \tensora \E, \A ). $ In \cite{majid_2} and \cite{beggsbrezinski}, this was done for bimodule connections on any differential calculus. Thus, if $g \in {\rm Hom}_{\A} ( \E \tensora \E, \A ) $ is a pseudo-Riemannian metric on $\E,$  we can make sense of the equality  $ \nabla_{( \E \tensora \E  )^*} g = 0 $ and we demonstrate that $ \nabla_{( \E \tensora \E  )^*} g = 0 $ if and only if $\nabla$ is compatible with $g$ in the sense of \cite{article3} ( and \cite{article1} ). In Section \ref{metriccompatibility}
, we recall a sufficient condition for the existence of a Levi-Civita connection on any tame differential calculus. This criterion will be crucially used in the next section.  In Section \ref{conformalsection}, we study Levi-Civita connections on conformal deformations of bilinear pseudo-Riemannian metrics.    In the last section, we define and compute the Ricci and scalar curvatures for some examples including an arbitrary `conformal deformation' of the canonical metric on the noncommutative $2$-torus as well as a bilinear metric on quantum Heisenberg manifolds.

\section{Some preliminaries on tame differential calculus} \label{formulation}

We begin by setting up the notations and conventions that are going to be followed. We will always work over the complex field. Unless otherwise mentioned, the symbol $\A$ will stand for a complex unital algebra and $ \Acenter $ will denote its center. The tensor product over the complex numbers $ \IC $ will be denoted by $ \tensorc $ while the notation $\tensora$ will denote the tensor product over the algebra $ \A. $ If $T$ is a linear map between suitable modules over $\A,$ $ {\rm Ran} ( T ) $ will denote the Range of $ T. $ 

The following well-known lemma will be useful for us.
\blmma \label{21stsep20}
If $\E$ is a finitely generated projective right $\A$-module, then there exist $ f^1, \cdots f^n $ in $\E$ and $ f_1, \cdots f_n $ in ${\rm Hom}_{\A} ( \E, \A ) $ such that for all $f \in \E$ and for all $ \phi \in {\rm Hom}_{\A} ( \E, \A ), $
	$$ f = \sum_i f^i f_i ( f ), ~ \phi = \sum_i \phi ( f^i ) f_i. $$
\elmma
$ \{ f^i, f_i: i = 1, \cdots n \} $ is known as a pair of dual bases for $\E.$  

Now we recall the definition of centered bimodules.
 \bdfn \label{centered}
We will say that a subset $ S  $ of a right $ \A $-module $ \E $ is right $ \A $-total in $ \E $ if   the right $\A$-linear span of $ S  $ equals $ \mathcal{E}. $ The center of  an $\A $-bimodule $ \mathcal{E} $ is defined to be the set $ \mathcal{Z} ( \mathcal{E} ) = \{ e \in \mathcal{E}: e a = a e ~ \forall ~ a ~ \in \A  \}. $ It is easy to see that $ \mathcal{Z} ( \E ) $ is a $ \mathcal{Z} ( \A ) $-bimodule.  $ \mathcal{E} $ is called centered if $ \mathcal{Z} ( \mathcal{E} ) $ is right $ \A $-total in $ \E. $
\edfn

If $\E$ and $\F$ are right $\A$-modules, $ {\rm Hom}_\A ( \mathcal{E}, \mathcal{F} ) $ will denote the set of all right $ \A $-linear maps from $ \E $ to $ \F. $ The set $ {\rm Hom}_\A ( \mathcal{E}, \mathcal{F} ) $ has a natural $\A$-bimodule structure which is as follows: 

the left $ \A $-module structure is given by left multiplication by elements of $ \A, $ i.e, for elements $ a $  in $ \A, ~ e  $ in $ \E $ and $ T  $ in $ {\rm Hom}_\A ( \mathcal{E}, \mathcal{F} ),$
 \begin{equation} \label{21stfeb20202} ( a. T ) ( e ) := a. T ( e ) \in \F. \end{equation} 
The right $ \A $ module structure on  $ {\rm Hom}_\A ( \mathcal{E}, \mathcal{F} ) $ is given by 
\begin{equation} \label{21stfeb20203} T. a ( e ) = T ( a e ). \end{equation}  
We will often make use of the following shorthand notation:      
\bdfn \label{Estar}
If $ \E $ is an $\A$-bimodule, then $ \E^* $ will stand for the $\A$-bimodule $ {\rm Hom}_{\A} ( \mathcal{E}, \A ).$ 
 \edfn 
The following isomorphism is well-known and will be used in the sequel:
\bdfn \label{xi}
Suppose $\E$ and $\F$ are finitely generated projective right $\A$-modules. The map $ \zeta_{\E, \F} $ will denote the canonical right $\A$-module isomorphism from $ \E \tensora \F^*  $ to $ \Hom_{\A} ( \F, \E ) $ which is defined by the following formula: 
$$ \zeta_{\E, \F} ( \sum_i e_i \tensora \phi_i ) ( f ) = \sum_i e_i \phi_i ( f ). $$
\edfn

Suppose $\A$ is an algebra. Then a differential calculus over $\A$ is a triplet $ ( \Omega^ ( \A ), \wedge, d ) $ where $ \Omega ( \A ) $ is a direct sum of $\A$-bimodules $ \Omega^j ( \A ) $, with $ \Omega^0 ( \A ) = \A$.    The map $ \wedge :  \Omega ( \A ) \tensora \Omega ( \A ) \rightarrow \Omega ( \A ) $ is an $\A$-bimodule map such that 
$$ \wedge ( \Omega^j ( \A ) \tensora \Omega^k ( \A )  ) \subseteq \Omega^{j + k} ( \A ).$$  
 $d$ is a map from $ \Omega^j ( \A ) $ to $ \Omega^{j + 1} ( \A ) $ such that 
				$$ d^2 = 0 ~ {\rm and} ~ d ( \omega \wedge \eta ) = d \omega \wedge \eta + ( - 1 )^{{\rm deg} ( \omega ) } \omega \wedge d \eta. $$
				Moreover, we will also assume that  $ \Omega^j ( \A ) $ is the right $ \A $-linear span of elements of the form $ d a_1 \wedge d a_1 \wedge \cdots \wedge d a_{j}. $
				
					We will often denote the bimodule of one-forms $\Omega^1 ( \A ) $ of a generic differential calculus by the symbol $\E.$ We will always assume that $\E$ is finitely generated and projective as a right $\A$-module. For notational convenience, we will sometimes denote a differential calculus by a pair $ ( \E, d ) $ if $\E$ is the bimodule of one-forms of a differential calculus $ ( \Omega ( \A ), \wedge, d ). $	
					
In Subsection \ref{tameproperties}, we define the notion of tame differential calculus and discuss some of its properties. In Subsection \ref{metricsontame}, we recall some results on pseudo-Riemannian metrics on tame differential calculi which will be useful for us.

\subsection{Definition and properties of tame differential calculi} \label{tameproperties}

		\bdfn \label{tame}
Suppose $\E$ is the bimodule of one-forms of a differential calculus $ ( \Omega ( \A ), \wedge, d ).$ We say that the differential calculus is tame if the following conditions hold:
\begin{enumerate}
    \item[i.] The bimodule $\E$ is finitely generated and projective as a right $\A$ module.
  	\item[ii.] The following short exact sequence of right $\A$-modules splits:
		$$ 0 \rightarrow {\rm Ker} ( \wedge ) \rightarrow \E \tensora \E \rightarrow \Omega^2 ( \A ) \rightarrow 0. $$
		Thus, in particular, there exists a right $\A$-module $\F$ isomorphic to $\twoform$ such that:
  \begin{equation} \label{splitting25thmay2018}
  \E \tensora \E = {\rm Ker}(\wedge) \oplus \F 
  \end{equation}
  \item[iii.] The map $u^\E:\Ecenter \otimes_{\Acenter} \A \rightarrow \E$ defined by
		$$u^\E(\sum_i e_i \otimes_{\Acenter} a_i)=\sum_i e_i a_i$$
		is an isomorphism of vector spaces,
		\item[iv.] Let us denote the idempotent in $ {\rm Hom}_{\A}(\E \tensora \E, \E \tensora \E)$ with range  ${\rm Ker}(\wedge)$ and  kernel $\F$ by the symbol $ \Psym $ and define $ \sigma = 2 \Psym - 1. $
We assume that $\sigma $ satisfies the following equation for all $\omega, \eta \in \Ecenter:$
		 \begin{equation} \label{17thdec20191} \sigma ( \omega \tensora \eta ) = \eta \tensora \omega. \end{equation}
	\end{enumerate}
\edfn
Let us note the following remark which will be very useful later on.
\brmrk \label{17thdec2019remark}
Proposition 2.4 of \cite{article1} asserts that the map $u^{\E} $ is actually a right $\A$-linear isomorphism from $\Ecenter \otimes_{\Acenter}\A $ to $\E.$ Moreover, in Lemma 4.4 of \cite{article1}, it has been proved that the map $ \Psym $ is $\A$-bilinear. Thus, the same is true about the map $\sigma = 2 \Psym - 1.$
\ermrk

Examples of tame differential calculi have been discussed in the introduction. We refer to Example 2.10 of \cite{article2} for more details.  In Section \ref{curvaturesection}, we will compute the scalar curvatures of the Levi-Civita connections for some tame differential calculi on the noncommutative torus and the quantum Heisenberg manifold.

	Let us recall some properties of a tame differential calculus from \cite{article1}, \cite{article3} and \cite{article2}. To begin with, Proposition 2.4 of \cite{article1} states that if $\E$ is the bimodule of one-forms of a tame differential calculus, then $\E$ is centered ( see Definition \ref{centered} ). Thus, as proved in  Section 4 of \cite{article3} ( see equation (2) ):
\begin{equation} \label{18thdec20191} a. e = e. a ~ {\rm if} ~ e \in \E ~ {\rm and} ~ a \in \Acenter.  \end{equation}
This implies that the bimodule of one-forms of a tame differential calculus is central in the sense of \cite{dubois} and \cite{dubois2}. Moreover, all properties of centered bimodules naturally continue to hold for a tame differential calculus. In particular, we have the following lemma.

\blmma \label{centeredremark} ( Lemma 1.4, \cite{article2} )
If $\E$ is the bimodule of one-forms of a tame differential calculus, then the following statements hold:
  \begin{enumerate}
 \item[i.]   $ \Ecenter $ is also left $\A$-total in $\E.$

\item[ii.]  The set $ \{ \omega \tensora \eta: \omega, \eta \in \Ecenter \} $ is both left and right $\A$-total in $\E \tensora \E.$

\item[iii.]  If $X$ is an element of $\E \tensora \E,$ there exist $ v_i $ in $\E,$ $ w_i \in \Ecenter $ and $ a_i $ in $\A$ such that 
   $$ X = \sum_i v_i \tensora w_i a_i. $$
\item[iv.]  If in addition, $\E$ is a free right $\A$-module with a basis $ \{ e_1, e_2, \cdots, e_n \} \subseteq \Ecenter, $ then any element $X$ in $\E$ can be written as a unique linear combination $ \sum_{i,j} e_i \tensora e_j a_{ij} $ for some elements $ a_{ij} $ in $\A.$
\end{enumerate} 	
\elmma
Now we explain the significance of the maps $\sigma$ and $\Psym.$ The map $\sigma$ plays the role of the flip map. In fact,
for all $\omega \in \Ecenter$ and $e \in \E,$ we have
		    \begin{equation} \label{10thjuly20182} \sigma ( \omega \tensora e ) = e \tensora \omega, ~ \sigma ( e \tensora \omega ) = \omega \tensora e \end{equation}
	and hence 
				\begin{equation} \label{21stmarch20} \Psym ( e \tensora \omega ) = \Psym ( \omega \tensora e ) = \frac{1}{2} ( \omega \tensora e + e \tensora \omega )\end{equation}
				 for all $\omega \in \Ecenter$ and $e \in \E.$
The decomposition $ \E \tensora \E = {\rm Ker} ( \wedge ) \oplus \mathcal{F} $ on simple tensors is explicitly given by
			\begin{equation} \label{21stmarch202} \omega \tensora \eta a = \frac{1}{2} ( \omega \tensora \eta a + \eta \tensora \omega a  ) + \frac{1}{2} ( \omega \tensora \eta a - \eta \tensora \omega a  ) \end{equation}
			for all $\omega, \eta$ in $\Ecenter$ and for all $a$ in $\A.$
			For the proof of these facts, we refer to Lemma 2.11 of \cite{article2}.

\subsection{Pseudo-Riemannian metrics on tame differential calculi} \label{metricsontame}

As discussed above, the map $\sigma = 2 \Psym - 1$ is an analogue of a flip map. Using this, we can now define pseudo-Riemannian metrics on a tame differential calculus. We will need the notation $ \E^*:= \Hom_\A ( \E, \A ) $ introduced in Definition \ref{Estar}.
\bdfn \label{metricdefn} ( \cite{article1}, \cite{article3} )
Suppose $\E$ is the bimodule of one-forms of a tame differential calculus $ ( \Omega ( \A ), \wedge, d ).$ A pseudo-Riemannian metric $ g $ on $ \E $ is
 an element of $ {\rm Hom}_{\A} ( \E \tensora \E, \A ) $ such that
 \begin{enumerate}
 \item[i.] $g$ is symmetric, i.e. $ g \circ \sigma = g, $
  \item[ii.] $g$  is non-degenerate, i.e,  the right $ \A$-linear map $ V_g: \E \rightarrow {\E}^* $ defined by $ V_g ( \omega ) ( \eta ) = g ( \omega \otimes_{\A} \eta ) $ is
 an isomorphism of right $ \A$-modules.
\end{enumerate}
We will say that a pseudo-Riemannian metric $g$ is a pseudo-Riemannian {\bf bilinear} metric if in addition, $ g $ is also left $\A$-linear.
\edfn
 Here, the right $ \A $-module structure on $ \E^* = \Hom_\A ( \E, \A ) $ is as in \eqref{21stfeb20203}, i.e, if $ \phi $ belongs to  $ \E^*, $ then for all $ a $ in  $ \A $ and $ e  $ in $\E,$ $ ( \phi. a ) ( e ) = \phi ( a. e ).  $

We collect some results about pseudo-Riemannian metrics on tame differential calculus. These had already been proven in some form or the other in \cite{article3} and \cite{article4}. We also refer to Lemma 2.12 of \cite{article2} for a complete proof.
\blmma \label{lemma0} 
 Suppose $ ( \E, d ) $ is a tame differential calculus and $g$ is a pseudo-Riemannian metric on $\E.$ Then the following are true:
\begin{itemize}
\item[i.] If  either $ e $  or $ f $ belongs to $ \Ecenter,$ we have
\begin{equation} \label{gsigmaisg} g (  e \tensora f ) = g ( f \tensora e ). \end{equation}
\item[ii.] If $ g_0 $ is a  pseudo-Riemannian  bilinear metric, then $ g_0 ( \omega \tensora \eta ) \in \mathcal{Z} ( \A ) $ if $ \omega, \eta  $ belong to $ \mathcal{Z} ( \E ).$
\item[iii.] If $e$ is an element of $\E$ such that $ g ( e \tensora \omega ) = 0 $ for all $ \omega  $ in $\Ecenter,$ then $e = 0.$ The same conclusion holds if $ g ( \omega \tensora e ) = 0 $ for all $ \omega $ in $\Ecenter.$
\end{itemize}
\elmma

\section{Connections and their compatibility with pseudo-Riemannian metrics} \label{nablag}

In this section, we study two definitions of metric-compatibility of connections on a tame differential calculus. In differential geometry, a connection $\nabla : \Omega^1 ( M ) \rightarrow \Omega^1 ( M ) \tensora \Omega^1 ( M ) $ canonically induces a connection $ \nabla_{( \Omega^1 ( M ) \otimes_{C^\infty ( M )} \Omega^1 ( M ) )^*} $ on $ ( \Omega^1 ( M ) \otimes_{C^\infty ( M )} \Omega^1 ( M ) )^*:= {\rm Hom}_{C^\infty ( M )} (  \Omega^1 ( M ) \otimes_{C^\infty ( M )} \Omega^1 ( M ), C^\infty ( M ) ).$ This is done by first extending $ \nabla $ to a connection on the bimodule $ \Omega^1 ( M ) \otimes_{C^\infty ( M )} \Omega^1 ( M ) $ and then using this connection to  define another connection on $ ( \Omega^1 ( M ) \otimes_{C^\infty ( M )} \Omega^1 ( M ) )^*.$

Unfortunately, when $\A$ is a noncommutative algebra and $\nabla$ is a connection on an $\A$-bimodule $\F,$ then there is no canonical construction of a connection on $\F \tensora \F.$ However, if $\nabla$ is a bimodule connection, then this can be done. We refer to Proposition 2.3 of \cite{majid_2} for a proof. Consequently, Beggs and Majid defined the compatibility condition of a connection $\nabla$ on $\F$ with the metric $g$ via the equality $ \nabla_{( \F \tensora \F )^*} g = 0. $ We refer to the paper \cite{majid_2} and also the book \cite{beggsmajidbook} for many interesting applications. 

We will show that if $ ( \E, d ) $ is a tame differential calculus, the $\A$-bimodule map $\sigma: \E \tensora \E \rightarrow \E \tensora \E $  introduced in Definition \ref{tame} helps us to circumvent the above mentioned problem. Thus, starting from an ( ordinary ) connection on $\E,$ there is a recipe to define a connection on $\E \tensora \E.$ This allows us to consider the equality $ \nabla_{( \F \tensora \F )^*} g = 0. $ We then show (  Theorem \ref{22ndsep203} ) that this equality holds if and only if the connection $\nabla$ on $\E$ is compatible with $g$ in the sense of \cite{article3} ( and \cite{article1} ).
 
We start with the definition of a connection and its torsion.
\bdfn ( \cite{frolich}, \cite{connes} ) \label{connection}
Suppose $ ( \Omega ( \A ), \wedge, d ) $ be a differential calculus on $\A$ and $\E:= \Omega^1 ( \A ). $   A (right) connection on an $\A$-bimodule $\F$
is a ${\mathbb C}$-linear map 
$$\nabla : \F   \rightarrow \F   \tensora \E  ~ {\rm satisfying} ~ {\rm the} ~ {\rm equation} ~  \nabla( f a) = \nabla( f )a + f \tensora da$$
for all $f$ in $\F   $ and $ a $ in $\A.$
 
The torsion of a connection $ \nabla $ on the bimodule $\E:= \Omega^1 ( \A ) $ is the right $\A$-linear map $ T_\nabla:= \wedge \nabla + d : \Omega^1 ( \A ) \rightarrow \Omega^2 ( \A ).$ The connection $\nabla$ is called torsionless if $ T_\nabla = 0. $
\edfn
For us, a connection will always mean a right connection, unless otherwise mentioned.
Let us recall that a tame differential calculus always admits a torsionless connection on the bimodule of one-forms.
\bthm  \label{torsionless} ( Theorem 3.3 of \cite{article3} )
Suppose $ ( \Omega ( \A ), \wedge, d ) $ is a tame differential calculus. Then the bimodule of one-forms $ \E = \Omega^1 ( \A ) $ admits a torsionless connection.
\ethm

In what follows, for an $\A$-bimodule $\F,$ the set $ \F^*:= {\rm Hom}_{\A} ( \F, \A )$ will be equipped with the $\A$-bimodule structure dictated by \eqref{21stfeb20202} and \eqref{21stfeb20203}.

We recall ( \cite{beggsbrezinski} ) that if $\nabla$ is a ( right ) connection on a finitely generated and projective right $\A$-module $\F,$ then it induces a left connection $\nabla_{\F^*}$ on $\F^*,$ i.e, $ \nabla_{\F^*} $ is a $\mathbb{C}$-linear map from $\F^*$ to $ \E \tensora \F^* $ such that for all $\phi$ in $\F^*$ and for all $a$ in $\A,$
$$ \nabla_{\F^*} ( a \phi ) = a \nabla_{\F^*} ( \phi ) + da \tensora \phi. $$
\bdfn ( Subsection 3.2, \cite{beggsbrezinski} ) \label{11thjan21}
If $\nabla: \F \rightarrow \F \tensora \E $ is a right connection on a finitely generated and projective right $\A$-module $\F$ with a pair of dual bases $ \{ f^i, f_i: i = 1,2, \cdots n \} $ as in Lemma \ref{21stsep20} and $ {\rm ev}: \E^* \tensora \E \rightarrow \A $ is the $\A$-bilinear map defined by
$ {\rm ev} ( \phi \tensora e ) = \phi ( e ), $
 then we define
$$ \nabla_{\F^*}: \F^* \rightarrow \E \tensora \F^*, ~ \nabla_{\F^*} ( \phi ) = \sum_i [ d ( \phi ( f^i ) ) \tensora f_i - ( {\rm ev} \tensora {\rm id} \tensora {\rm id} ) ( \phi \tensora \nabla ( f^i ) \tensora f_i )  ].$$
\edfn  

The definition of $ \nabla_{\F^*} $ is independent of the choice of $ \{ f^i, f_i: i = 1,2, \cdots n \} $ as can be seen from the next proposition.

	\bppsn \label{19thsep20} ( \cite{beggsbrezinski} ) 
Suppose $ ( \E, d ) $ is a any differential calculus and $\nabla$ a connection on a finite generated and projective right $\A$-module $\F.$ If $ \zeta_{\E, \F}: \E \tensora \F^* \rightarrow {\rm Hom}_{\A} ( \F, \E ) $ is the isomorphism as in Definition \ref{xi} and ${\rm ev}: \E^* \tensora \E \rightarrow \A $ the $\A$-bilinear map as above, then for all $\phi \in \F^*$ and all $f \in \E,$ we get 
\begin{equation} \label{22ndsep20} \zeta_{\E, \F} ( \nabla_{\F^*} ( \phi ) ) ( f ) = d ( \phi ( f ) ) - ( \ev \tensora \id ) ( \phi \tensora \nabla ( f ) ). \end{equation}
Thus, the definition of $\nabla_{\F^*}$ is independent of the choice of the dual bases. Moreover, $\nabla_{\F^*}$ is a left connection on $\F^*.$ 
\eppsn
{\bf Proof:} This follows from the proof of Proposition 2.9 of \cite{majid_2}. However, for the sake of completeness, we give a proof. Indeed, for $\phi$ in $\F^*$ and $f \in \E,$ we make the following computation:
\begin{eqnarray*}
&& \zeta_{\E, \F} ( \nabla_{\F^*} ( \phi ) ) ( f )\\
&=& ( \id \tensora \ev ) ( \nabla_{\F^*} ( \phi ) \tensora f )\\
&=& \sum_i [ d ( \phi ( f^i ) ) f_i ( f ) - ( \ev \tensora \id ) ( \phi \tensora \nabla ( f^i ) f_i ( f ) )]\\
&=& \sum_i [  d ( \phi ( f^i ) f_i ( f ) ) - \phi ( f^i ) d ( f_i ( f ) ) - ( \ev \tensora \id ) ( \phi \tensora \nabla ( f^i f_i ( f ) ) ) + ( \ev \tensora \id ) ( \phi \tensora f^i d ( f_i ( f ) ) )   ]\\
&=&  d ( \phi ( f ) ) - \sum_i \phi ( f^i ) d ( f_i ( f ) )  - ( \ev \tensora \id ) ( \phi \tensora \nabla ( f ) ) + \sum_i \phi ( f^i ) d ( f_i ( f ) )\\
&=& d ( \phi ( f ) ) - ( \ev \tensora \id ) ( \phi \tensora \nabla ( f ) ).
\end{eqnarray*}
This proves the first assertion. Next, we prove that $\nabla_{\F^*}$ is a left connection on $\F^*:$
\begin{eqnarray*}
&& \nabla_{\F^*} ( a \phi )\\
&=& \sum_i [ d ( a \phi ( f^i ) ) \tensora f_i - ( {\rm ev} \tensora {\rm id} \tensora {\rm id} ) ( a \phi \tensora \nabla ( f^i ) \tensora f_i )  ] \\
&=&  \sum_i a [ d ( \phi ( f^i ) ) \tensora f_i - ( {\rm ev} \tensora {\rm id} \tensora {\rm id} ) ( \phi \tensora \nabla ( f^i ) \tensora f_i )  ] + \sum_i da \phi ( f^i ) \tensora f_i\\
&& {\rm (} ~ {\rm by} ~ {\rm the} ~ {\rm Leibniz} ~ {\rm rule} ~ {\rm and} ~ {\rm the} ~ {\rm left} ~ \A-{\rm linearity} ~ {\rm of} ~ {\rm the} ~ {\rm map} ~ {\rm ev} ~ {\rm )}\\
&=& a \nabla_{\F^*} ( f ) + da \tensora \phi
\end{eqnarray*} 
by the second equation of Lemma \ref{21stsep20}. This completes the proof. 
\qed

Now we define the notion of the compatibility of a connection with a pseudo-Riemannian metric on a tame differential calculus. We will need the following proposition:
	
\bppsn ( Subsection 4.1, \cite{article1} ) \label{30thmarch203}
If $g$ is a pseudo-Riemannian metric on the bimodule of one-forms $\E$ of a tame differential calculus $ ( \E, d ), $ we define $\Pi_g^0(\nabla):\Ecenter\tensorc\Ecenter \rightarrow \E$ as the map given by
$$\Pi_g^0(\nabla)(\omega \tensorc \eta)=(g\tensora {\rm id})\sigma_{23}(\nabla(\omega)\tensora \eta + \nabla(\eta) \tensora \omega).$$
Then $ \Pi_g^0 $ extends to a well defined map from $ \E \tensora \E $ to $ \E $ to be denoted by $ \Pi_g ( \nabla ).$ 
\eppsn
It turns out that for $\omega, \eta $ in $\Ecenter$ and $ a $ in $\A,$ the following equation holds: 
 \begin{equation} \label{20thfeb2020} \Pi_g ( \nabla ) ( \omega \tensora \eta a ) = \Pi^0_g ( \nabla ) ( \omega \otimes_{\Acenter} \eta ) a + g ( \omega \tensora \eta ) da. \end{equation}
We recall that  by Lemma \ref{centeredremark}, any element of $\E \tensora \E $ is a finite sum of elements of the form $ \omega \tensora \eta a, $ where $\omega, \eta \in \Ecenter$ and $a$ belongs to $\A.$  Therefore, the equation \eqref{20thfeb2020} defines the map $ \Pi_g ( \nabla ) $ on the whole of $ \E \tensora \E.$ 

Now we are in a position to define the metric-compatibility of a connection in our set-up.
\bdfn \label{19thfeb202021}
Suppose $ ( \E, d ) $ is a tame differential calculus and $ g $ is a pseudo-Riemannian metric on $\E.$

 A connection $\nabla$ on $\E$ is said to be compatible with $g$ if 
$$\Pi_g(\nabla) ( e \tensora f ) = d ( g ( e \tensora f ) )~ {\rm for} ~ {\rm all} ~ e,~ f ~ {\rm in} ~ \E.$$
A connection $\nabla$ on $\E$ which is torsionless and compatible with $g$ is called a Levi-Civita connection for the triplet $(\E, d, g ).$
\edfn

Let us introduce a Sweedler-type notation. If $ \nabla $ is a connection on $\E$ and $ e  $ belongs to $\E,$ then we will write
\begin{equation} \label{19thsep202} \nabla ( e ) = e_\zero \tensora e_\one.  \end{equation}
So if $\omega, \eta \in \Ecenter$ and $a \in \A,$ then
$$ \Pi_g ( \nabla ) ( \omega \tensora \eta a ) = g ( \omega_\zero \tensora \eta ) \omega_\one a + g ( \eta_\zero \tensora \omega ) \eta_\one a + g ( \omega \tensora \eta ) da.$$
By ii. of Lemma \ref{centeredremark}, we know that any element of $\E \tensora \E$ is a finite linear combination of terms of the form $ \omega \tensora \eta a $ where $\omega, \eta \in \Ecenter $ and $a \in \A. $ Hence, a connection $\nabla$ is compatible with $g$ if and only if
\begin{equation} \label{19thsep203} d ( g ( \omega \tensora \eta a ) ) = g ( \omega_\zero \tensora \eta ) \omega_\one a + g ( \eta_\zero \tensora \omega ) \eta_\one a + g ( \omega \tensora \eta ) da  \end{equation}
 for all $\omega, \eta \in \Ecenter $ and for all $a \in \E. $

Now we demonstrate that if $\nabla$ is any ( right ) connection on the space of one-forms of a tame differential calculus $ ( \E, d ), $ then we can lift $\nabla$ to a connection on $\E \tensora \E.$ We do this in two steps: in the first step, we define a map
$$\widetilde{\nabla}_{\E \tensora \E}: ( \Ecenter \otimes_{\Acenter} \Ecenter ) \rightarrow ( \E \tensora \E ) \tensora \E. $$
In the second step, we define a map $\nabla_{\E \tensora \E}: ( \E \tensora \E ) \rightarrow ( \E \tensora \E ) \tensora \E.$

Indeed, we define
$$\widetilde{\nabla}_{\E \tensora \E}: \Ecenter \otimes_{\Acenter} \Ecenter \rightarrow \E \tensora \E, ~ \widetilde{\nabla}_{\E \tensora \E} ( \omega \otimes_{\Acenter} \eta ) = \sigma_{23} ( \nabla ( \omega ) \tensora \eta ) + \omega \tensora \nabla ( \eta ). $$
We need to check that $ \widetilde{\nabla}_{\E \tensora \E} $ is well-defined. Thus, for $\omega, \eta $ in $\Acenter$ and $a$ in $\Acenter,$ we compute
\begin{eqnarray*}
 \nabla_{\E \tensora \E} ( \omega a \otimes_{\Acenter} \eta ) &=& \sigma_{23} ( \nabla ( \omega a ) \tensora \eta ) + \omega a \tensora \nabla ( \eta )\\
&=& \sigma_{23} ( \nabla ( \omega ) \tensora a \eta ) + \omega \tensora \eta \tensora da + \omega \tensora a \nabla ( \eta ) ~ {\rm (} ~ {\rm by} ~  \eqref{10thjuly20182} ~ {\rm )}\\
&=& \sigma_{23} ( \nabla ( \omega ) \tensora a \eta ) + \omega \tensora \eta \tensora da + \omega \tensora  \nabla ( \eta ) a   ~ {\rm (} ~ {\rm by} ~ \eqref{18thdec20191} ~ {\rm )}\\
&=& \sigma_{23} ( \nabla ( \omega ) \tensora a \eta ) + \omega \tensora  \nabla ( \eta a )\\ 
&=& \sigma_{23} ( \nabla ( \omega ) \tensora a \eta ) + \omega \tensora  \nabla ( a \eta  ) ~ {\rm (} ~ {\rm as} ~ \eta \in \Ecenter ~ {\rm )}\\
&=& \nabla_{\E \tensora \E} ( \omega \otimes_{\Acenter} a \eta )
\end{eqnarray*}
proving that $ \widetilde{\nabla}_{\E \tensora \E} $ is well-defined.
Now we execute the second step. From Definition 4.7 of \cite{article1},  we know that the map 
$$ {u}^{\E \tensora \E}: \Ecenter \otimes_{\Acenter} \Ecenter \otimes_{\Acenter} \A \rightarrow \E \tensora \E,~ {u}^{\E \tensora \E} ( \omega \otimes_{\Acenter} \eta \otimes_{\Acenter} a ) = \omega \tensora \eta a $$
 is an isomorphism. So it makes sense to define
$$ \nabla_{\E \tensora \E}: ( \E \tensora \E ) \rightarrow ( \E \tensora \E ) \tensora \E ~ {\rm by} ~ {\rm the} ~ {\rm formula} $$ 
$$ \nabla_{\E \tensora \E} (  {u}^{\E \tensora \E} ( \omega \otimes_{\Acenter} \eta \otimes_{\Acenter} a )    ) = \widetilde{\nabla}_{\E \tensora \E} ( \omega \otimes_{\Acenter} \eta ) a + \omega \tensora \eta \tensora da. $$

\bppsn \label{22ndsep202}
If $ ( \E, d ) $ is a tame differential calculus and $\nabla$ a connection on $\E,$ then  $ \nabla_{\E \tensora \E} $ is a well-defined connection on  $ \E \tensora \E. $
\eppsn
{\bf Proof:} We start by proving that $\nabla_{\E \tensora \E}$ is well-defined.  From the defining formula of $ \nabla_{\E \tensora \E}, $ it is clear that it suffices to prove
\begin{equation} \label{20thsep20} \nabla_{\E \tensora \E} ( \omega \otimes_{\Acenter} \eta b \otimes_{\Acenter} a ) = \nabla_{\E \tensora \E} ( \omega \otimes_{\Acenter} \eta \otimes_{\Acenter} b a )  \end{equation}
for all $\omega, \eta \in \Ecenter, $ $a \in \A$ and for all $b \in \Acenter.$ We compute
\begin{eqnarray*}
\nabla_{\E \tensora \E} ( {u}^{\E \tensora \E}  ( \omega \otimes_{\Acenter} \eta b \otimes_{\Acenter} a ) )  &=&  \widetilde{\nabla}_{\E \tensora \E} ( \omega \tensora \eta b ) a + \omega \tensora \eta \tensora b da \\
&=& \sigma_{23} ( \nabla ( \omega ) \tensora \eta b ) a + \omega \tensora \nabla ( \eta b ) a + \omega \tensora \eta \tensora b da \\
&=& \widetilde{\nabla}_{\E \tensora \E} ( \omega \otimes_{\Acenter} \eta ) b a + \omega \tensora \eta \tensora d ( b a ) \\
&& ~ {\rm (}  ~ {\rm we} ~ {\rm have} ~ {\rm applied} ~ {\rm the} ~ {\rm Leibniz} ~ {\rm rules} ~ {\rm for} ~ \nabla ~ {\rm and} ~ d ~ {\rm )}\\ 
&=& \nabla_{\E \tensora \E} ( {u}^{\E \tensora \E}  ( \omega \otimes_{\Acenter} \eta \otimes_{\Acenter} b a ) ).
\end{eqnarray*}

Now we prove that $\nabla_{\E \tensora \E}$ is a connection. 
If $\omega, \eta \in \Ecenter $ and $a, b \in \A, $ we obtain
\begin{eqnarray*}
\nabla_{\E \tensora \E} ( ( \omega \tensora \eta a ) b ) &=& \nabla_{\E \tensora \E} {u}^{\E \tensora \E} ( \omega \otimes_{\Acenter} \eta \otimes_{\Acenter} a b )\\
&=& \widetilde{\nabla}_{\E \tensora \E} ( \omega \otimes_{\Acenter} \eta ) ab + \omega \tensora \eta \tensora d ( a b )\\
&=&  \widetilde{\nabla}_{\E \tensora \E} ( \omega \otimes_{\Acenter} \eta ) ab + \omega \tensora \eta \tensora da b + \omega \tensora \eta \tensora a db\\
&=& \nabla_{\E \tensora \E} ( \omega \tensora \eta a ) b + ( \omega \tensora \eta a ) \tensora db.
\end{eqnarray*}
This completes the proof.
\qed

Summarizing, if $ ( \E, d ) $ is a tame differential calculus and $\nabla$ a ( right ) connection on $\E,$ then by Proposition \ref{22ndsep202}, we have a ( right ) connection $ \nabla_{\E \tensora \E} $ on $\E \tensora \E.$ Then Proposition \ref{19thsep20} delivers a left connection $\nabla_{( \E \tensora \E )^*} $ on $ ( \E \tensora \E )^*.$ If $g \in \Hom_{\A} ( \E \tensora \E, \A ) $ is a pseudo-Riemannian metric on $\E,$ we can therefore 
 make sense of the quantity $\nabla_{( \E \tensora \E )^*}.$

If $ \{ f^i, f_i : i = 1,2, \cdots n  \} $ is a pair of dual bases for a finitely generated and projective right $\A$-module $\F,$ then it can be easily checked that $ \{ f^i \tensora f^j, f_i \tensora f_j: i,j = 1,2, \cdots n \} $ is a pair of dual bases for $\F \tensora \F.$ Therefore, Definition \ref{11thjan21} implies that
$$ \nabla_{( \E \tensora \E )^*} g = \sum_{i,j} [  d (  g ( f^i \tensora f^j ) ) \tensora ( f_i \tensora f_j ) - ( \ev \tensora \id \tensora \id ) ( g \tensora \nabla_{\E \tensora \E} ( f^i \tensora f^j ) \tensora ( f_i \tensora f_j ) )  ].$$

Classically, it is well-known that a connection $\nabla$ on the cotangent bundle $ \Omega^1 ( M ) $ of a manifold $M$ is compatible with a pseudo-Riemannian metric $g$ if and only if $ \nabla_{( \Omega^1 ( M ) \otimes_{C^\infty ( M )} \Omega^1 ( M ) )^*} g = 0. $ We show that the same is true with our definition of metric-compatibility as in Definition \ref{19thfeb202021}. 

\bthm \label{22ndsep203}
if $ ( \E, d ) $ is a tame differential calculus and $g$ is a pseudo-Riemannian metric on $\E,$ then a connection $\nabla$ on $\E$ is compatible with $g$ in the sense of Definition \ref{19thfeb202021} if and only if $ \nabla_{( \E \tensora \E )^*} g = 0. $
\ethm
{\bf Proof:} The proof follows from  \eqref{22ndsep20}. Since $ \{ \omega \tensora \eta a: \omega, \eta \in \Ecenter, a \in \A \}  $ is right $\A$-total in $\E \tensora \E $ ( part ii. of Lemma \ref{centeredremark} ) and $ \zeta_{\E, \E \tensora \E}: \E \tensora ( \E \tensora \E  )^* \rightarrow {\rm Hom}_{\A} ( \E \tensora \E, \E ) $ is an isomorphism ( Definition \ref{xi} ),  $ \nabla_{( \E \tensora \E )^*} g = 0 $ if and only if 
$$ \zeta_{\E, \E \tensora \E} ( \nabla_{( \E \tensora \E )^*}  g ) ( \omega \tensora \eta a  )  = 0$$
for all $\omega, \eta \in \Ecenter $ and for all $a \in \A.$
The equation \eqref{22ndsep20} implies that $ \nabla_{( \E \tensora \E )^*} g = 0 $ if and only if for all $\omega, \eta$ and $a$ as above,
$$ d ( g ( \omega \tensora \eta a ) ) = ( \ev \tensora \id ) ( g \tensora \nabla_{\E \tensora \E} ( \omega \tensora \eta a ) ).$$
But 
\begin{eqnarray*}
&& ( \ev \tensora \id ) ( g \tensora \nabla_{\E \tensora \E} ( \omega \tensora \eta a ) )\\
&=& ( \ev \tensora \id ) ( g \tensora \nabla_{\E \tensora \E} ( \omega \tensora \eta ) a + g \tensora \omega \tensora \eta \tensora da )\\
&=& ( \ev \tensora \id ) ( g \tensora \sigma_{23} ( \nabla ( \omega ) \tensora \eta ) a + g \tensora \omega \tensora \nabla ( \eta ) a ) + g ( \omega \tensora \eta ) da
\end{eqnarray*}
by the definition of $ \nabla_{\E \tensora \E}. $ Using the Sweedler-type notation introduced in \eqref{19thsep202}, the above expression is equal to 
\begin{eqnarray*}
 && g ( \omega_\zero \tensora \eta ) \omega_\one a + g ( \omega \tensora \eta_\zero ) \eta_\one a + g ( \omega \tensora \eta ) da\\
 &=&  g ( \omega_\zero \tensora \eta ) \omega_\one a + g ( \eta_\zero \tensora \omega ) \eta_\one a + g ( \omega \tensora \eta ) da,
\end{eqnarray*}
since $\omega \in \Ecenter$ and we have applied \eqref{gsigmaisg}.
Therefore, $ \nabla_{( \E \tensora \E )^*} g = 0 $ if and only if for all $\omega, \eta \in \Ecenter $ and for all $a \in \A, $
$$ d ( g ( \omega \tensora \eta a ) ) = g ( \omega_\zero \tensora \eta ) \omega_\one a + g ( \eta_\zero \tensora \omega ) \eta_\one a + g ( \omega \tensora \eta ) da. $$
Comparing with \eqref{19thsep203}, we deduce that $\nabla$ is compatible with $g$ if and only if $ \nabla_{( \E \tensora \E  )^*} g = 0. $
\qed

\section{A criterion for the existence of Levi-Civita connections} \label{metriccompatibility}

In \cite{article3} ( also see \cite{article1} ) and \cite{article2}, existence and uniqueness of Levi-Civita connections have been proved for for bilinear and strongly $\sigma$-compatible pseudo-Riemannian  metrics  respectively. We will take the path adopted in \cite{article1} and \cite{article2} to study Levi-Civita connections for conformally deformed pseudo-Riemannian metrics in Section \ref{conformalsection}. This will require a sufficient condition for the existence of Levi-Civita connections ( Theorem \ref{19thfeb2020} ). 

We will need the following definition: 
\bdfn \label{19thfeb20202}
For a tame differential calculus $ ( \E, d ), $ the symbol $ \E \tensorsym \E$ will denote   ${\rm Ker} ( \wedge ) = {\rm Ran} ( \Psym ). $ If $g$ is a pseudo-Riemannian metric, the element $dg$ will denote the map
$$ dg: \E \tensora \E \rightarrow \E, dg ( e \tensora f ) = d ( g ( e \tensora f )  ).$$
\edfn
Then we have the following theorem:
\bthm ( Theorem 4.14, \cite{article1} ) \label{19thfeb2020}
 Let $ ( \E, d ) $ be a tame differential calculus and $g$ a pseudo-Riemannian metric on $\E.$ Let $ \E \tensorsym \E $ be the $\A$-bimodule of Definition \ref{19thfeb20202}.  
$$ {\rm We} ~ {\rm define} ~ {\rm a} ~ {\rm map} ~ \Phi_g : {\rm Hom}_{\A}(\E, \E \tensorsym \E) \rightarrow {\rm Hom}_{\A}(\E \tensorsym \E, \E) ~ {\rm by} ~ {\rm the} ~ {\rm formula:} $$
$$ \Phi_g(L) ( X ) =(g \tensora {\rm id}) \sigma_{23} (L \tensora {\rm id})(1+\sigma) ( X ) $$
for all $X$ in $\E \tensorsym \E.$

Then $ \Phi_g $ is right $\A$-linear. Moreover, if $ \Phi_g:{\rm Hom}_{\A}(\E, \E \tensorsym \E) \rightarrow {\rm Hom}_{\A}(\E \tensorsym \E, \E) $ is an isomorphism of right $\A$-modules, then there exists a unique connection $ \nabla $ on $\E$ which is torsion-less and compatible with $g.$
Moreover, if $\nabla_0$ is a fixed torsionless connection on $\E,$ then  $ \nabla $ is given by the following equation: 
				\be \label{26thnov2} \nabla = \nabla_0 + \Phi^{- 1}_g ( dg - \Pi_g ( \nabla_0 ) ).  \ee
	Here, $dg: \E \tensora \E \rightarrow \E $ is the map defined in Definition \ref{19thfeb20202}.			
\ethm

The proof of this theorem works for any pseudo-Riemannian metric. The formula \eqref{26thnov2} follows from the proof of Theorem 4.14 of \cite{article1}. We only need to remark that the proof of Theorem 4.13 of \cite{article1} uses the existence of a torsion-less connection on $\E.$ In our case, this condition is satisfied by virtue of Theorem \ref{torsionless}.

The main result of \cite{article3} and \cite{article1} is the following:

\bthm  ( Theorem 6.1 of \cite{article3}, Theorem 4.1 of \cite{article1} ) \label{21stfeb2020}
Let $ ( \E, d ) $ be a tame differential calculus and $g_0$ be a pseudo-Riemannian {\bf bilinear} metric on $\E.$ Then there exists a unique Levi-Civita connection for the triplet $(\E, d, g_0 ).$
\ethm
\brmrk \label{22ndfeb2020}
In \cite{article1}, Theorem \ref{21stfeb2020} was proved by verifying that the map $ \Phi_{g_0}: {\rm Hom}_{\A}(\E, \E \tensorsym \E) \rightarrow {\rm Hom}_{\A}(\E \tensorsym \E, \E) $ is an isomorphism of right $\A$-modules.
\ermrk
 In \cite{article3}, a completely different proof was given. Indeed, the uniqueness of such a connection followed by deriving a Koszul type formula of a torsionless and $g_0$ compatible connection. The existence followed by proving that the above mentioned Koszul-formula indeed defines a torsionless and $g_0$-compatible connection on $\E.$
 We have been unable to generalize the proof of \cite{article3} for  metrics which are not $\A$-bilinear. However, we demonstrate that for conformal deformations of a {\bf bilinear} pseudo-Riemannian metrics on a tame differential calculus, Theorem \ref{21stfeb2020} allows us to give a short proof for existence and uniqueness of Levi-Civita connection. This is the content of the next section. In \cite{article2}, Theorem \ref{21stfeb2020} was generalized to the case of strongly $\sigma$-compatible pseudo-Riemannian metrics.
 
\section{The Levi-Civita connection for a conformally deformed metric} \label{conformalsection}
		 
		Let $(\E, d )$ be a tame differential calculus and $g_0$ be a pseudo-Riemannian {\bf bilinear} metric on $\E.$ We fix an invertible element  $k$ of $\A.$
		\bdfn \label{7thjune20}
		With $ (\E, d, g_0 ) $ as above, the map 
		$$g: \E \tensora \E \rightarrow \A,~ g ( e \tensora f ) = k g ( e \tensora f )  $$
		is called a conformal deformation of $g_0.$
		\edfn
		Indeed, it can be easily checked that $g$ is a pseudo-Riemannian metric on $\E.$ 
		
		Throughout this section, we will follow the notations developed till now so that $g$ will denote a conformal deformation of the pseudo-Riemannian {\bf bilinear} metric $g_0$ on a tame differential calculus $ ( \E, d ).$ Moreover, the map $dg: \E \tensora \E \rightarrow \E $ will be as defined in Definition \ref{19thfeb20202}. We note that if $ g = k. g_0 $ is a conformal deformation of $g_0,$ then for all $e, f$ in $\E,$
		$$ dg ( e \tensora f ) = d ( g ( e \tensora f ) ) = d ( k.g_0 ( e \tensora f  ) ) = dk. g_0 ( e \tensora f ) + k. d g_0 ( e \tensora f ) $$
		and so 
		\begin{equation} \label{18thmarch20}  dg = dk. g_0 + k. d g_0 ~ {\rm as} ~ {\rm maps} ~ {\rm from} ~ \E \tensora \E ~ {\rm to} ~ \E. \end{equation}
		This section has three results. Firstly, in Theorem \ref{19thjan1}, we prove that there exists a unique Levi-Civita connections for the triplet $ (\E, d, g ).$ Secondly, in Theorem \ref{19thjan4}, we derive a concrete formula for this Levi-Civita connection in terms of $k$ and a fixed torsionless connection on $\E.$ Finally, Proposition \ref{20thjanprop} deduces the Christoffel symbols for the Levi-Civita connection for the conformal deformation $g$ when $\E$ is in addition a free right $\A$-module satisfying some conditions.

	Let us clarify a couple of notations to be used in Theorem \ref{19thjan1}.
  We recall  that $ \E \tensorsym \E $ is defined to be $ {\rm Ker} ( \wedge ). $ However, from the definition of a differential calculus, we know that the map $ \wedge: \E \tensora \E \rightarrow \twoform $ is an $\A$-bimodule map and so $ \E \tensorsym \E =  {\rm Ker} ( \wedge ) $ is an $\A$-bimodule.
 Therefore, from \eqref{21stfeb20202}, we know that $ \Hom_{\A} ( \E \tensorsym \E, \E ) $ is a left ( as well as a right ) $\A$-module. For an invertible element $k$ in $\A,$ let				
					$$ L_k: {\rm Hom}_\A ( \E \tensorsym \E, \E ) \rightarrow {\rm Hom}_\A ( \E \tensorsym \E, \E ), ~ L_k ( X ) = k.X $$
					denote the left $\A$-module multiplication on $ {\rm Hom}_\A ( \E \tensorsym \E, \E ).$
					
For the pseudo-Riemannian bilinear metric $g_0$ as above, 
$$ dk. g_0: \E \tensora \E \rightarrow \E ~ {\rm will} ~ {\rm denote} ~ {\rm the} ~ {\rm map} ~ {\rm defined} ~ {\rm by} ~ dk. g_0 ( e \tensora f ) = dk. g_0 ( e \tensora f ). $$
Since $g_0$ is right $\A$-linear, it is clear that $ dk. g_0 $ is an element of $ \Hom_{\A} ( \E \tensora \E, \E ). $ Since $ \E \tensorsym \E = {\rm Ker} ( \wedge ) $ is both a left and right $\A$-submodule of $\E \tensora \E$ by the above discussions, it makes sense to view the element $ dk. g_0 $ as an element of $ \Hom_{\A} ( \E \tensorsym \E, \E ).$ Thus, $ L_{k^{-1}} (  dk. g_0 ) $ is also an element of $ \Hom_{\A} ( \E \tensorsym \E, \E ).$ Therefore, by Remark \ref{22ndfeb2020}, we can conclude that $ \Phi^{- 1}_{g_0} L_{k^{-1}} ( dk. g_0 ) $ is a well-defined element of $ \Hom_{\A} ( \E, \E \tensorsym \E ). $ We will use these facts in Theorem \ref{19thjan1}.					
					
	We have the following result:
					\bthm \label{19thjan1} 
Let $ ( \E, d ) $ be a tame differential calculus and $g_0$ a pseudo-Riemannian bilinear metric on $\E.$ We will denote the Levi-Civita connection for the triplet $ ( \E, d, g_0  ) $ by the symbol $ \nabla_{g_0}. $ If $k$ is an invertible element of $\A,$ then there exists a unique Levi-Civita connection $\nabla$ for the triplet $ ( \E, d, k. g_0 )$ given by the formula:					
					\begin{equation} \label{22ndfeb20202} \nabla = \nabla_{g_0} + \Phi^{- 1}_{g_0} L_{k^{-1}} ( dk. g_0 ). \end{equation}
					Here, $ \Phi_{g_0} $ is the map defined in Theorem \ref{19thfeb2020}.
					\ethm
					{\bf Proof:} 
					As stated above, $g$ will denote the pseudo-Riemannian metric $k. g_0.$ We use Theorem \ref{19thfeb2020} to prove the existence and uniqueness of Levi-Civita connection for the triplet $ ( \E, d, g ). $ Thus, it suffices to prove that the map
					$$ \Phi_g: \Hom_{\A} ( \E, \E \tensorsym \E ) \rightarrow \Hom_{\A} ( \E \tensorsym \E, \E ) $$
					is a right $\A$-linear isomorphism. But since $ g = k. g_0, $ it is easy to verify that $ \Phi_g = L_k. \Phi_{g_0}. $ By Remark \ref{22ndfeb2020}, $ \Phi_{g_0} $ is a right $\A$-linear isomorphism from $ \Hom_{\A} ( \E, \E \tensorsym \E ) $ to $ \Hom_{\A} ( \E \tensorsym \E, \E ) $ and so by Theorem \ref{19thfeb2020}, the Levi-Civita connection $\nabla_{g_0}$ exists.
					
					Since $ k $ is an invertible element of $\A,$ $ \Phi_g = L_k \Phi_{g_0} $ is also invertible and its inverse is explicitly given by
					\begin{equation} \label{22ndfeb20203} ( L_k. \Phi_{g_0} )^{-1} = \Phi^{-1}_{g_0} ( L_k )^{-1} = \Phi^{-1}_{g_0} L_{k^{-1}}. \end{equation}
					In particular, the hypothesis of Theorem \ref{19thfeb2020} is satisfied and we have a unique Levi-Civita connection for the triplet $ ( \E, d, g ).$
					
					Next, the equation \eqref{22ndfeb20202} follows from \eqref{26thnov2}. Indeed, the Levi-Civita connection $ \nabla_{g_0} $ for the triplet $ ( \E, d, g_0 ) $ is torsionless and so \eqref{22ndfeb20203} implies that
					\begin{equation} \label{22ndfeb20204} \nabla = \nabla_{g_0} + \Phi^{-1}_{g_0} L_{k^{-1}} ( dg - \Pi_g ( \nabla_{g_0} ).  \end{equation}
					Now since $ g = k. g_0, $ by using \eqref{20thfeb2020}, it can be easily checked that 
					$$ \Pi_g ( \nabla_{g_0} ) = k. \Pi_{g_0} ( \nabla_{g_0} ) = k. dg_{0} $$
					since $ \nabla_{g_0} $ is compatible with the metric $ g_0 $ ( Definition \ref{19thfeb202021} ). Hence, \eqref{22ndfeb20204} implies that 
					\begin{eqnarray} \nabla &=& \nabla_{g_0} + \Phi^{-1}_{g_0} L_{k^{-1}} ( dg - k. d g_0  ) \nonumber \\
					                        &=& \nabla_{g_0} + \Phi^{-1}_{g_0} L_{k^{-1}}  (  dk. g_0 + k. d g_0 - k. d g_0    ) ~ {\rm (} ~ {\rm by} ~ \eqref{18thmarch20} ~ {\rm )}\\
																	&=& \nabla_{g_0} 	+ \Phi^{-1}_{g_0} L_{k^{-1}}  (  dk. g_0 ) \nonumber
	                                				\end{eqnarray}					
				This finishes the proof of the theorem. 
					\qed
				
				The formula \eqref{22ndfeb20202} for the Levi-Civita connection $\nabla$ in Theorem \ref{19thjan1} can be made more explicit by deriving a formula for $ \Phi^{-1}_{g_0}. $ This is the content of Theorem \ref{19thjan4} for which we need to need to make a definition. 
				\bdfn		\label{omegag0defn}
			For a pseudo-Riemannian bilinear metric $ g_0, $ we  define an element $ \Omega_{g_0} \in \E \tensora \E $ by  
			$$ \Omega_{g_0} = ( {\rm id}_\E \tensora V^{- 1}_{g_0}  ) \zeta^{- 1}_{\E,\E} ( {\rm id}_\E ).  $$
				\edfn
	We recall that $ \zeta_{\E, \E}: \Hom_{\A} ( \E, \E ) \rightarrow \E \tensora \E^*  $ is the the right $\A$-module isomorphism defined in Definition \ref{xi}. Moreover, as $g_0$ is {\bf bilinear}, it can be easily checked that the map $ V_{g_0} $ ( and hence   $ V^{-1}_{g_0} $ ) is left $\A$-linear. Therefore, the map $  \id_{\E} \tensora V^{-1}_{g_0} $ makes sense.		
			
				Now we are ready to state the following Theorem: 
				\bthm \label{19thjan4}
				Suppose $ ( \E, d ) $ is a tame differential calculus and $g_0$ is a pseudo-Riemannian bilinear metric on $\E.$ If $ \nabla_{g_0} $ is the Levi-Civita connection for $ ( \E, d, g_0 ) $ and $k$ is an invertible element of $\A,$ then the Levi-Civita connection $\nabla$ of $ ( \E, d, k g_0 ) $ is given by the following formula:
$$ \nabla ( \omega ) = \nabla_{g_{0}}(\omega) + k^{-1} \Psym ( dk \tensora \omega ) - \frac{1}{2} k^{-1} \Omega_{g_0} g_0 ( dk \tensora \omega ).$$
Here, the map $ \Psym: \E \tensora \E \rightarrow \E \tensora \E $ is the one defined in Definition \ref{tame}.
\ethm
The proof of this theorem will be derived in steps in the next subsection. The Theorem \ref{19thjan4} will help us to derive the Christoffel symbol of the Levi-Civita connection when the module is free.

	\subsection{A formula for the inverse of $\Phi_{g_0}$}
	
	In this subsection, we prove Theorem \ref{19thjan4}. We will continue with the notations made before. In particular, we will be using the map $\Psym$ introduced in Definition \ref{tame} while $ \Omega_{g_0} $ is as defined in Definition \ref{omegag0defn}. 	
	
	Comparing the statements of Theorem \ref{19thjan1} and Theorem \ref{19thjan4}, it is clear that we need to prove the following equation for all $e$ in $\E:$
	$$ \Phi^{-1}_{g_0} ( L_{k^{-1}} ( dk. g_0 ) ) ( e ) = k^{-1} \Psym ( dk \tensora e ) - \frac{1}{2} k^{-1} \Omega_{g_0} g_0 ( dk \tensora e ). $$

\blmma \label{omegag0sym}
				$ \Omega_{g_0} $ is an element of  $ \E \tensorsym \E. $ 
				\elmma
				{\bf Proof:} We will need the  following $\A$-bilinear  map from \cite{article3} ( also see \cite{article1} ):
				\begin{equation} \label{22ndfeb207} g^{(2)}_0 : ( \E \tensora \E ) \tensora ( \E \tensora \E ) \rightarrow \cla,  ~ g^{(2)}_0 ( ( \omega \tensora \eta   ) \tensora ( \omega^\prime \tensora \eta^\prime ) ) = g_0 ( \omega g_0 ( \eta \tensora \omega^\prime ) \tensora \eta^\prime ).\end{equation}
				Then by Proposition 6.6 of \cite{article3} ( also see Proposition 3.8 of \cite{article1} ), we know that
		 \begin{equation} \label{14thdec20191} g^{(2)}_0 ( \theta \tensora \theta^{\prime} ) = 0 ~ \forall ~ \theta^\prime \in \E \tensora \E ~ {\rm implies} ~ {\rm that} ~  \theta = 0. \end{equation} 
		Moreover, by Lemma 4.17 of \cite{article1}, 
		\begin{equation} \label{14thdec20192} g^{(2)}_0 ( \sigma ( e \tensora f ) \tensora ( e^\prime \tensora f^\prime ) ) = g^{(2)}_0 ( ( e \tensora f ) \tensora \sigma ( e^\prime \tensora f^\prime  ) ) \end{equation}
		for all $e,f, e^\prime, f^\prime$ belonging to $\E.$		
				
				Since $ \E \tensorsym \E $ is the range of the idempotent $  \Psym  $ by definition, we need to prove that  $ \Psym ( \Omega_{g_0} ) = \Omega_{g_0}.$ Since $ \sigma = 2 \Psym - 1  $ ( Definition \ref{tame} ), this amounts to proving the equation 
			 $ \sigma \Omega_{g_0} = \Omega_{g_0}.$
			
	We claim 		that for all $ \omega^\prime, \eta^\prime \in \mathcal{Z} ( \E ), $
				\begin{equation} \label{14thdec20194} g^{(2)}_0 ( \sigma \Omega_{g_0} \tensora ( \omega^\prime  \tensora \eta^\prime  )   ) = g^{(2)}_0 (  \Omega_{g_0} \tensora ( \omega^\prime  \tensora \eta^\prime  )   ) .\end{equation}
				If  \eqref{14thdec20194} is true, then  by the right $\A$-linearity of $ g^{(2)}_0,$ we get that for all $ \omega^\prime, \eta^\prime \in \Ecenter $ and for all $a$ in $\A,$
				$$ g^{(2)}_0 ( ( \sigma \Omega_{g_0} - \Omega_{g_0} ) \tensora ( \omega^\prime \tensora \eta^\prime a )     ) = 0. $$
				By  part ii. of Lemma \ref{centeredremark}, we deduce that
				$$ g^{(2)}_0 ( ( \sigma \Omega_{g_0} - \Omega_{g_0}   ) \tensora \theta^\prime  ) = 0 ~ {\rm for} ~ {\rm all} ~ \theta^\prime ~ {\rm in} ~ \E \tensora \E. $$
		Therefore, \eqref{14thdec20191} implies that $ \sigma \Omega_{g_0} = \Omega_{g_0}. $ Thus, we are left to prove \eqref{14thdec20194}.
		
		Once again we use Lemma \ref{centeredremark} to recall that there exist $ v_i $ in $\E,$ $ w_i  $ in $ \Ecenter $ and $ a_i $ in $\A$ such that
		\begin{equation} \label{22ndfeb206} \Omega_{g_0} = \sum_i v_i \tensora w_i a_i.   \end{equation}
		On the other hand, by the definition of $ \Omega_{g_0}, $ for all $e$ in $\E,$ we obtain 
		\begin{eqnarray}
		 e &=& ( \zeta_{\E, \E} ( {\rm id}_{\E} \tensora V_{g_0}  ) \Omega_{g_0}   ) ( e ) \nonumber \\
		   &=& ( \zeta_{\E, \E} ( {\rm id}_{\E} \tensora V_{g_0}  ) (  \sum_i v_i \tensora w_i a_i   )   ) ( e ) \nonumber \\
			 &=& \sum_i \zeta_{\E, \E} (  v_i \tensora V_{g_0} (  w_i a_i  )   ) ( e )  \nonumber \\
			 &=& \sum_i v_i V_{g_0} ( w_i a_i ) ( e ) ~ {\rm (} ~ {\rm by} ~ {\rm Definition} ~ \ref{xi} ~ {\rm )} \nonumber \\
		   &=& \sum_i v_i g_0 ( w_i a_i \tensora e ). \label{omegag01} 
			\end{eqnarray}
			If $ \omega^\prime, \eta^\prime \in \mathcal{Z} ( \E ), $ then
				\begin{eqnarray*}
				 g^{(2)}_0 ( \sigma \Omega_{g_0} \tensora ( \omega^\prime \tensora \eta^\prime ) ) &=&  g^{(2)}_0 ( \Omega_{g_0} \tensora \sigma ( \omega^\prime \tensora \eta^\prime  ) ) ~ {\rm (} ~ {\rm by} ~ \eqref{14thdec20192} ~ {\rm )}\\
				&=&  g^{(2)}_0 (  \Omega_{g_0} \tensora ( \eta^\prime \tensora \omega^\prime ) )  ~ {\rm (} ~ {\rm by} ~ \eqref{17thdec20191} ~ {\rm )}\\ 
&=& 	g^{(2)}_0 ( \sum_i ( v_i \tensora w_i  a_i ) \tensora  ( \eta^\prime \tensora \omega^\prime  ) )~ {\rm (} ~ {\rm by} ~ \eqref{22ndfeb206} ~ {\rm )}\\
&=& g_0 (  \sum_i v_i g_0 ( w_i a_i \tensora \eta^\prime  ) \tensora \omega^\prime   )~ {\rm (} ~ {\rm by} ~ \eqref{22ndfeb207} ~ {\rm )}\\
&=& g_0 ( \eta^\prime \tensora \omega^\prime ) ~ {\rm (} ~ {\rm by} ~ \eqref{omegag01} ~ {\rm )}\\
&=& g_0 \circ \sigma ( \eta^\prime \tensora \omega^\prime )~ {\rm (} ~ {\rm by} ~ {\rm Definition} ~ \ref{metricdefn} ~ {\rm )}\\
&=& g_0 ( \omega^\prime \tensora \eta^\prime  ) ~ {\rm (} ~ {\rm by} ~ \eqref{17thdec20191} ~ {\rm )}\\
&=& g_0 ( \sum_i v_i g_0 ( w_i a_i \tensora \omega^\prime  ) \tensora \eta^\prime ) ~ {\rm (} ~ {\rm by} ~ \eqref{omegag01} ~ {\rm )}\\
&=& g^{(2)}_0 ( \Omega_{g_0} \tensora ( \omega^\prime \tensora \eta^\prime ) )~ {\rm (} ~ {\rm by} ~ \eqref{22ndfeb207} ~ {\rm )}.			
				\end{eqnarray*}
				 This finishes the proof of \eqref{14thdec20194} and hence the lemma.
			 	\qed

				\blmma \label{28thnov2}
				For all $\eta$ in $\Ecenter,$ the following equation holds:
				$$ ( g_0 \tensora {\rm id} ) \sigma_{23} ( \Omega_{g_0} \tensora \eta ) = \eta. $$
				\elmma
				{\bf Proof:} Let us continue writing $ \Omega_{g_0} $ as $ \sum_i v_i \tensora w_i a_i $ ( finitely many terms ) for some $ v_i \in \E, w_i \in \clz ( \E ) $ and $ a_i \in \A $ as in Lemma \ref{omegag0sym} so that the relation $ \sigma \Omega_{g_0} = \Omega_{g_0} $ ( as obtained from Lemma \ref{omegag0sym} ) implies that
				\be \label{28thnov3} \sum_i v_i \tensora w_i a_i = \sigma ( \sum_i v_i \tensora \omega_i a_i   ) = \sum_i w_i \tensora v_i a_i \ee
			as $\sigma$ is right $\A$-linear and we have applied the second equation of \eqref{10thjuly20182}.
			Moreover, if $ \eta $ in $ \mathcal{Z} ( \E ), $ we have
			  \begin{eqnarray}
			  ( g_0 \tensora {\rm id} ) \sigma_{23} ( \Omega_{g_0} \tensora \eta ) &=& ( g_0 \tensora {\rm id} ) \sigma_{23} ( \sum_i v_i \tensora w_i a_i \tensora \eta   )  \nonumber \\
				&& ( g_0 \tensora {\rm id} ) ( \sum_i v_i \tensora \eta \tensora w_i a_i ) ~ {\rm (} ~ {\rm by} ~ \eqref{10thjuly20182}     {\rm )}  \nonumber \\ 
				&=& \sum_i g_0 ( v_i \tensora \eta ) w_i a_i. \label{15thdec20193}
				\end{eqnarray}
				  If $ \eta^\prime $ belongs to $\Ecenter,$ we compute:
				\begin{eqnarray*}
				g_0(  ( g_0 \tensora {\rm id} ) \sigma_{23} ( \Omega_{g_0} \tensora \eta )  \tensora \eta^\prime ) &=& g_0 ( \sum_i g_0 ( v_i \tensora \eta ) w_i a_i  \tensora \eta^\prime )~ {\rm (} ~ {\rm by} ~ \eqref{15thdec20193} ~ {\rm )}\\
				 &=& g_0 ( \sum_i g_0 ( v_i \tensora \eta ) w_i   \tensora a_i \eta^\prime ) \\
				&=&g_0 ( \sum_i w_i g_0 (  v_i \tensora \eta )    \tensora a_i \eta^\prime )~ {\rm (} ~ {\rm as} ~ \omega_i \in \Ecenter ~ {\rm )}\\
				&=& g^{(2)}_0 ( \sum_i ( w_i \tensora v_i  ) \tensora ( \eta \tensora a_i \eta^\prime  )  )~ {\rm (} ~ {\rm by} ~ \eqref{22ndfeb207} ~ {\rm )}\\
				&=& g^{(2)}_0 ( \sum_i ( w_i \tensora v_i  ) a_i \tensora ( \eta \tensora  \eta^\prime  )  )~ {\rm (} ~ {\rm as} ~ \eta \in \Ecenter ~ {\rm )}\\
				&=&  g^{(2)}_0 ( \sum_i ( v_i \tensora w_i a_i  ) \tensora ( \eta \tensora \eta^\prime  )  ) ~ ( {\rm by} ~ \eqref{28thnov3} )\\
				&=& g_0 ( \sum_i v_i g_0 ( w_i a_i \tensora \eta ) \tensora \eta^\prime )\\
				&=& g_0 ( \eta \tensora \eta^\prime ) ~ ( {\rm by} ~ \eqref{omegag01} ).
				\end{eqnarray*}
				Therefore, for all $\eta^\prime$ in $\Ecenter,$ we obtain
				$$ g_0 ( (  ( g_0 \tensora \id  ) \sigma_{23} ( \Omega_{g_0} \tensora \eta ) - \eta  )   \tensora \eta^\prime  ) = 0. $$
				By part iii of Lemma \ref{lemma0}, we can conclude that $ ( g_0 \tensora {\rm id} ) \sigma_{23} ( \Omega_{g_0} \tensora \eta ) = \eta. $ This proves the lemma.
				\qed

				Having obtained the above results, we are now in a position to prove Theorem \ref{19thjan4}.

{\bf Proof of Theorem \ref{19thjan4}:} Let us define the map
						$$T_{dk}: \E \tensora \E \rightarrow \E ~ {\rm by} ~ T_{dk} ( \omega \tensora \eta ) = dk g_0 ( \omega \tensora \eta ). $$
					By an abuse of notation, we will denote the restriction of $ T_{dk}  $ to $ \E \tensorsym \E $ by the same symbol $ T_{dk}. $ We claim that the following equation holds:
						\begin{equation} \label{19thmarch20} ( \Phi^{-1}_{g_0} ( T_{dk} )  ) ( \omega ) = \Psym ( dk \tensora \omega ) - \frac{1}{2} \Omega_{g_0} g_0 ( dk \tensora \omega ).\end{equation} 
	If \eqref{19thmarch20} holds, then the theorem follows from a computation.  Indeed, we observe that the map $ \Phi_{g_0} $ is left $\A$-linear since $ g_0 $ is so and therefore,
\begin{eqnarray*}
 &&\Phi^{-1}_{g_0} ( k^{-1} dk. g_0 ) ( \omega )\\
 &=& k^{-1} \Phi^{-1}_{g_0} ( dk. g_0 ) ( \omega )\\
 &=& k^{-1} (  \Psym ( dk \tensora \omega   ) - \frac{1}{2} \Omega_{g_0} g_0 (  dk \tensora \omega   )    ) ~ {\rm (} ~  {\rm by} ~ \eqref{19thmarch20}  ~ {\rm )} 
\end{eqnarray*}
so that by Theorem \ref{19thjan1},
 $$ \nabla ( \omega ) = \nabla_{g_0} ( \omega ) + k^{-1} \Psym ( dk \tensora \omega ) - \frac{1}{2} k^{-1} \Omega_{g_0} g_0 ( dk \tensora \omega ). $$
Thus, we are left to prove our claim, i.e, \eqref{19thmarch20}.

 We define $L^{1}_{dk}\in{\rm Hom}_{\cla}(\cle, \cle\ot_{\cla}^{{\rm sym}}\cle)$ as
		\begin{eqnarray*}
		 L^{1}_{dk}(\omega)= \Psym ( dk \ot_{\cla}\omega)-\frac{1}{2}\Omega_{g_{0}}g_{0}( dk \ot_{\cla}\omega).
		\end{eqnarray*}
		Since $ \Omega_{g_0} $ belongs to $ \E \tensorsym \E $ by Lemma \ref{omegag0sym}, $ L^1_{dk} ( \omega ) $ indeed belongs to $ \E \tensorsym \E = {\rm Ran} ( \Psym ). $
		
                We want to prove $\Phi_{g_{0}}L^{1}_{dk}=T_{dk}$. So for $\omega,\eta\in\clz(\cle)$, we compute 
                \begin{eqnarray*}
                 && \Phi_{g_{0}}L^{1}_{dk}(\omega\ot_{\cla}\eta)\\
								 &=& (g_{0}\ot_{\cla}{\rm id})\sigma_{23}(L^{1}_{dk}\ot_{\cla}{\rm id})(\omega\ot_{\cla}\eta+\eta\ot_{\cla}\omega)\\
                 &=&(g_{0}\ot_{\cla}{\rm id})\sigma_{23}(L^{1}_{dk}(\omega)\ot_{\cla}\eta+L^{1}_{dk}(\eta)\ot_{\cla}\omega)\\
                 &=& (g_{0}\ot_{\cla}{\rm id})\sigma_{23}( \Psym (dk \tensora \omega)\ot_{\cla}\eta-\frac{1}{2}\Omega_{g_{0}}g_{0}(dk \ot_{\cla}\omega)\ot_{\cla}\eta + \Psym (dk \ot_{\cla}\eta) \tensora \omega\\
								&-& \frac{1}{2}\Omega_{g_{0}}g_{0}(dk \ot_{\cla}\eta)\ot_{\cla}\omega).
								\end{eqnarray*}
								Since $\omega, \eta$ belong to $\Ecenter,$ we can apply \eqref{21stmarch20}  to rewrite the above expression as:
	\begin{eqnarray*}
  && \frac{1}{2}(g_{0}\ot_{\cla}{\rm id})\sigma_{23}(dk \ot_{\cla}\omega\ot_{\cla}\eta+\omega\ot_{\cla}dk \ot_{\cla}\eta + dk \ot_{\cla}\eta\ot_{\cla}\omega+\eta\ot_{\cla} dk \ot_{\cla}\omega-\Omega_{g_{0}}g_{0}(dk \ot_{\cla}\omega)\ot_{\cla}\eta\\
	&-& \Omega_{g_{0}}g_{0}(dk \ot_{\cla}\eta)\ot_{\cla}\omega)\\
  &=& \frac{1}{2}[g_{0}(dk \ot_{\cla}\eta)\omega+g_{0}(\omega\ot_{\cla}\eta)dk + g_{0}( dk \ot_{\cla}\omega)\eta + g_{0}(\eta\ot_{\cla}\omega)dk - (g_{0}\ot_{\cla}{\rm id})\sigma_{23}(\Omega_{g_{0}}\ot_{\cla}\eta)g_{0}(dk \ot_{\cla}\omega)\\
	&-& (g_{0}\ot_{\cla}{\rm id})\sigma_{23}(\Omega_{g_{0}}\ot_{\cla}\omega)g_{0}(dk \ot_{\cla}\eta)]
  \end{eqnarray*}
	and in the last step we have used \eqref{10thjuly20182} as well as the fact that $\omega, \eta$ belong to $\Ecenter.$
	
Now using Lemma \ref{28thnov2}, \eqref{gsigmaisg} and the fact that $\omega,\eta\in\clz(\cle)$, the expression reduces to $2g_{0}(\omega \ot_{\cla} \eta)dk. $ However, since $g_{0}$ is a bilinear metric and $\omega,\eta\in\clz(\cle)$, by part ii. of Lemma \ref{lemma0} and \eqref{18thdec20191}, 
$$g_{0}(\eta\ot_{\cla}\omega) dk= dk g_{0}(\omega \ot_{\cla} \eta).$$
 Hence for all $\omega,\eta\in\clz(\cle)$,
$$\Phi_{g_{0}}( L^{1}_{dk} ) (\omega \tensora \eta)=T_{dk}(\omega \tensora \eta).$$
Since the set $ \{  \omega \tensora \eta: \omega, \eta \in \Ecenter    \} $ is right $\cla$-total in $\cle\ot_{\cla}\cle$ ( part ii. of Lemma \ref{centeredremark} ) and the maps $ \Phi_{g_0} ( L^1_{dk} ), T_{dk} $ are right $\A$-linear, we can finally conclude that for all $e, f $ in $\E,$
$$ \Phi_{g_0} ( L^1_{dk} ) ( e \tensora f ) = T_{dk} ( e \tensora f ).$$
This proves our claim and hence the theorem.
 \qed

\subsection{Christoffel symbols for a class of free modules}

We end this section by Proposition \ref{20thjanprop} which computes the Christoffel symbols of the Levi-Civita connection for a diagonal metric on a class of free modules. The hypotheses of Proposition \ref{20thjanprop} ( and Proposition  \ref{20thmarch203} ) are satisfied for the differential calculi on the noncommutative torus and the quantum Heisenberg manifold  for which  we refer to Subsection \ref{subnctorus} and Subsection \ref{heisenbergcurvature}.  On the way, we will prove Proposition \ref{20thmarch203} which establishes a sufficient condition on the calculus $ ( \E, d ) $ such that the Christoffel symbols are symmetric ( see \eqref{20thjan2} ).
We start by recalling the definition of the Christoffel symbols of a connection.
					
					\bdfn	 \label{christoffel}
						Suppose that $ \E $ is a free module with a basis $ \{ e_1, e_2, \cdots e_n \} \in \mathcal{Z} ( \E ) $ and $ \nabla $ is a connection on $\E.$ Then we can define the ``Christoffel symbols" $  \Gamma^i_{jk} \in \cla $ as follows:
					\be \label{19thjan2} \nabla ( e_i ) = \sum_{j,k} e_j \tensora e_k \Gamma^i_{jk}. \ee
					\edfn
	We note that since $ \nabla ( e_i ) $ belongs to $\E \tensora \E, $ the elements $ \Gamma^i_{jk} $ are uniquely defined by part iv. of Lemma \ref{centeredremark}.

					\bppsn \label{20thmarch203}
					If  $ ( \E, d ) $ is a tame differential calculus such that $\E$ is a free right $\A$-module with a basis $ \{ e_1, e_2, \cdots e_n \} \in \mathcal{Z} ( \E ) $ such that   $ d ( e_i ) = 0 $ for all $ i = 1,2, \cdots n.$ Then we have the following:
				\begin{enumerate}
				\item[i.]  There exist derivations $ \partial_j : \E \rightarrow \cla, ~ j = 1,2, \cdots n, $ such that
					\be \label{19thjan3} d a = \sum_j e_j \partial_j ( a ). \ee
				\item[ii.]  The Christoffel symbols of a torsion-less connection satisfy
					\be \label{20thjan2}  \Gamma^p_{kl} = \Gamma^p_{lk} ~ {\rm for} ~ {\rm all} ~ p,k,l. \ee 
					\end{enumerate}
					\eppsn
				{\bf Proof:} Let $ a $ be in element of $\A.$ Hence, $ da $ belongs to $\E.$ Since $ ( \E, d ) $ is a tame differential calculus, $\E$ is centered ( see Subsection \ref{tameproperties} ) and  so there exist unique elements $ a_i $ in $\A$ such that
				$$ da = \sum_j e_j a_j. $$
	For all $j = 1,2, \cdots n,$ we define 
	$$ \partial_j ( a ):= a_j.$$
We need to check that $ \partial_j $ is a derivation for all $j.$ So we fix two elements $a$ and $b$ in $\A.$ Then by the definition of $ \partial_j,$ we have
\begin{eqnarray*}
\sum_j e_j \partial_j ( a.b ) &=& d ( a.b ) =  da. b + a. db ~ {\rm (}~ {\rm since} ~ d ~ {\rm is} ~ {\rm a} ~ {\rm derivation} ~ {\rm )}\\
                            &=& \sum_j e_j \partial_j ( a ). b + a \sum_j  e_j \partial_j ( b ) = \sum_j e_j \partial_j ( a ). b + \sum_j e_j a. \partial_j ( b )~ {\rm (}~ {\rm since} ~ e_j ~ {\rm belongs} ~ {\rm to} ~ \Ecenter ~ {\rm )}\\ 
													&=& \sum_j e_j ( \partial_j ( a ). b + a. \partial_j ( b ) )	
\end{eqnarray*}
By comparing the coefficients of $e_j,$ we conclude that for all $j = 1,2, \cdots n,$
$$ \partial_j ( a. b ) = \partial_j ( a ). b + a. \partial_j ( b ), $$
i.e, $\partial_j$	is a derivation.

Now we prove the second assertion. We observe that by \eqref{21stmarch20},
$$ \Psym ( e_i \tensora e_j + e_j \tensora e_i ) = \frac{1}{2} ( e_i \tensora e_j + e_j \tensora e_i ) + \frac{1}{2} ( e_j \tensora e_i + e_i \tensora e_j   ) = e_i \tensora e_j + e_j \tensora e_i.$$
Therefore, $ e_i \tensora e_j + e_j \tensora e_i  $ belongs to $ {\rm Ran} ( \Psym ) = {\rm Ker} ( \wedge ) $ ( by Definition \ref{tame} ).	Hence, 
\begin{equation} \label{20thmarch204} e_i \wedge e_j = \wedge ( e_i \tensora e_j ) = - \wedge ( e_j \tensora e_i ) = - e_j \wedge e_i. \end{equation}
Since $ \nabla $ is assumed to be torsionless and $ d ( e_i ) = 0, $ we get
\begin{eqnarray*}
0 &=& - d ( e_i ) = \wedge \circ \nabla ( e_i )\\
  &=& \wedge  ( \sum_{j,k} e_j \tensora e_k \Gamma^i_{jk}  ) = \sum_{j,k} e_j \wedge e_k \Gamma^i_{jk}\\
	&=& \sum_{j < k} e_j \wedge e_k ( \Gamma^i_{jk} - \Gamma^i_{kj}   )
\end{eqnarray*} 				
since $e_i \wedge e_i = 0 $ by \eqref{20thmarch204}. Therefore, 
$$ \sum_{j < k} e_j \wedge e_k ( \Gamma^i_{jk} - \Gamma^i_{kj} ) = 0. $$
Thus, the proof will be complete once we prove that $ \{ e_j \wedge e_k: j < k \} $	is an $\A$-linearly independent set in $ \Omega^2 ( \A ).$

Consider the set $ \{ \frac{1}{2} ( e_i \tensora e_j - e_j \tensora e_i  ): i < j  \} \subseteq \E \tensora \E$ which is linearly independent ( over $\A$ ) as $ \{ e_i \tensora e_j: i,j  \} $ is a basis of the free right $\A$-module $\E \tensora \E.$ By \eqref{21stmarch202}, the set $ \{ \frac{1}{2} ( e_i \tensora e_j - e_j \tensora e_i  ): i < j  \} $ is actually contained in the right $\A$-module $ \mathcal{F}.$ Since $ ( \E, d ) $ is tame, the splitting of the short exact sequence in part ii. of Definition \ref{tame} implies that the map $ \wedge: \F \rightarrow \Omega^2 ( \A ) $ is a right $\A$-module isomorphism.

Hence, $ \wedge ( \{ \frac{1}{2} ( e_i \tensora e_j - e_j \tensora e_i   ): i < j  \} ) $ is an $\A$-linearly independent set in $ \Omega^2 ( \A ). $ But if $i < j,$
$$ \wedge (  \frac{1}{2} ( e_i \tensora e_j - e_j \tensora e_i   )  ) = e_i \wedge e_j $$ 
by \eqref{20thmarch204}. This proves that $ \{ e_j \wedge e_k: j < k \} $ is an $\A$-linearly independent set in $\Omega^2 ( \A ) $ and completes the proof of the proposition.     			 
				\qed

	\brmrk
	The condition $ d ( e_i ) = 0 $ is necessary for the equation \eqref{20thjan2} to hold. Indeed, consider the computation of the Levi-Civita connection for the fuzzy sphere in Section 8 of \cite{article3}. From Remark 8.7 of that paper, it is evident that the Christoffel symbols do not satisfy the relation \eqref{20thjan2} while equation ( 36 ) of \cite{article3} shows that $ d ( e_m ) \neq 0. $
				
This is a completely noncommutative phenomenon since in classical differential geometry, the equation \eqref{20thjan2} is always satisfied. Indeed, in the classical case, the Christoffel symbols are defined on a local chart $ ( U, x ) $ and the cotangent bundle is free over the open set $U$ with a basis $ \{ e_i:= d x_i: i = 1, \cdots, n \},$ $n$ being the dimension of the manifold. Hence, $ d e_i = d^2 x_i =0. $ 				\ermrk

					\bppsn \label{20thjanprop}
				Suppose  $ ( \E, d ) $ is a tame differential calculus and $g_0$ is a bilinear pseudo-Riemannian metric such that the following conditions are satisfied:
					\begin{enumerate}
				\item $ \E $ is a free right $\A$-module with a basis $ \{ e_1, e_2, \cdots e_n \} \in \mathcal{Z} ( \E ) $ such that   $ d ( e_i ) = 0 $ for all $ i = 1,2, \cdots n.$
			
				\item $ g_0 ( e_i \tensora e_j ) = \delta_{ij} $ for all $i,j.$
				\end{enumerate}
						
		We will denote the Christoffel symbols of the Levi-Civita connection $\nabla_{g_0}$ for the triplet $ ( \E, d, g_0 ) $ by the symbol $ ( \Gamma_0 )^i_{jl}.$ Consider the conformally deformed metric $ g:= k. g_0 $ where $k$ is an invertible element in $\cla.$ Then the  Christoffel symbols of the Levi-Civita connection $\nabla$  for the triplet $ ( \E, d, g ) $   are given by: 
		\be \label{22ndjan2} \Gamma^i_{jl} = ( \Gamma_0 )^i_{jl} + \frac{1}{2} ( \delta_{il} k^{- 1} \partial_j ( k ) +  \delta_{ij} k^{- 1} \partial_l ( k ) -  \delta_{jl} k^{- 1} \partial_i ( k )  ). \ee
		Here, $ \partial_i $ are the derivations as in \eqref{19thjan3}.
\eppsn
		{\bf Proof:} 
		To begin with, we note that the existences of $ \nabla_{g_0} $ and $ \nabla_0 $ follow from Theorem \ref{21stfeb2020} and Theorem \ref{19thjan4} respectively.
			
			We claim that under our assumptions, 
		\begin{equation} \label{20thmarch20} \Omega_{g_0} = \sum_i e_i \tensora e_i. \end{equation}
		Indeed, part iv. of Lemma \ref{centeredremark} implies the existence of elements $a_{ij}$ in $\A$ such that 
		$$ \Omega_{g_0} = \sum_{i,j} e_i \tensora e_j a_{ij}. $$
		Now for a fixed  $ k_0, $ we have
		\begin{eqnarray*}
   e_{k_0}&=&  ( g_0 \tensora {\rm id} ) \sigma_{23} ( \Omega_{g_0} \tensora e_{k_0} ) ~ {\rm (} ~ {\rm by} ~ {\rm Lemma} ~ \ref{28thnov2} ~ {\rm )}\\
	    &=& \sum_{i,j} ( g_0 \tensora {\rm id} ) \sigma_{23} ( e_i \tensora e_j a_{ij} \tensora e_{k_0} )\\
	&=& ( g_0 \tensora {\rm id} ) ( \sum_{i,j} e_i \tensora e_{k_0} \tensora e_j a_{ij} ) ~ {\rm (} ~ {\rm by} ~ {\rm the} ~ {\rm second} ~ {\rm equation} ~ {\rm of} ~ \eqref{10thjuly20182} )\\
	&=& \sum_j e_j a_{k_0 j}
	\end{eqnarray*}
	as $ g_0 ( e_i \tensora e_{k_0} ) = \delta_{i, k_0}.$ Therefore, we can deduce that $ a_{k_0 k_0}= 1  $ and $ a_{k_0 j} = 0 $ if $ j \neq k_0. $ This proves the claim.
	
	Now we apply Theorem \ref{19thjan4} to see that
	\begin{eqnarray} 
	\nabla ( e_i ) &=& \nabla_{g_0} ( e_i ) + k^{-1} \Psym ( dk \tensora e_i ) - \frac{1}{2} k^{-1} \Omega_{g_0} g_0 ( dk \tensora e_i ) \nonumber \\ 
	 &=& \nabla_{g_0} ( e_i ) +  \frac{1}{2} k^{-1} dk \tensora e_i + \frac{1}{2} k^{-1} e_i \tensora dk - \frac{1}{2} k^{-1} \Omega_{g_0}  g_0 ( dk \tensora e_i ), \label{20thmarch202} 
	\end{eqnarray}
	where, we have applied \eqref{21stmarch20}. Now, since $e_i$ belongs to $\Ecenter,$ $k. e_i = e_i . k$ and so
	\begin{eqnarray*}
	k^{-1} \Omega_{g_0} &=& k^{-1} ( \sum_i e_i \tensora e_i ) ~ {\rm (} ~ {\rm by} ~ \eqref{20thmarch20} ~ {\rm )}\\
											&=& ( \sum_i e_i \tensora e_i ) k^{-1}\\
											&=& \Omega_{g_0} k^{-1}.
	\end{eqnarray*}
	Therefore,
	\begin{eqnarray*}
	&& \frac{1}{2} k^{-1} dk \tensora e_i + \frac{1}{2} k^{-1} e_i \tensora dk - \frac{1}{2} k^{-1} \Omega_{g_0}  g_0 ( dk \tensora e_i )\\
	&=& \frac{1}{2} k^{-1} dk \tensora e_i + \frac{1}{2} k^{-1} e_i \tensora dk - \frac{1}{2}  \Omega_{g_0} k^{-1}  g_0 ( dk \tensora e_i )\\
	&=& \frac{1}{2} k^{-1} dk \tensora e_i + \frac{1}{2} k^{-1} e_i \tensora dk - \frac{1}{2} \sum_l e_l \tensora e_l g_0 ( k^{-1} dk \tensora e_i   )\\
	&& {\rm (} ~ {\rm since} ~ g_0 ~ {\rm is} ~ {\rm left} ~ \A-{\rm linear} ~ {\rm and} ~ {\rm we} ~ {\rm have} ~ {\rm applied} ~ \eqref{20thmarch20} ~ {\rm )}\\
	&=& \frac{1}{2} ( \sum_j k^{- 1} e_j \partial_j ( k ) \tensora e_i + \sum_j e_i \tensora k^{ - 1 } e_j \partial_j ( k )  ) - \frac{1}{2} ( \sum_l e_l \tensora e_l g_0 ( k^{ -1} ( \sum_j e_j \partial_j ( k )  ) \tensora e_i  ) )\\
	  && {\rm (} ~ {\rm by} ~ \eqref{19thjan3} ~ {\rm and} ~ {\rm as} ~ e_i ~ \in  ~ \Ecenter ~ {\rm )}\\
	   &=& \sum_{j,l} e_j \tensora e_l ( \frac{1}{2} \delta_{il} k^{- 1} \partial_j ( k ) + \frac{1}{2} \delta_{ij} k^{- 1} \partial_l ( k ) -   \frac{1}{2} \delta_{jl} k^{- 1} \partial_i ( k )   )
		\end{eqnarray*}
		as $ g_0 ( e_j \tensora e_i ) = \delta_{ij}, e_i \in \Ecenter $ and $g_0$ is $\A$-bilinear.
			
		Hence, by \eqref{20thmarch202}, we obtain
		\begin{eqnarray*}
		 \nabla ( e_i ) &=&  \sum_{j,l} e_j \tensora e_l (  \Gamma_0  )^i_{jl}  + \sum_{j,l} e_j \tensora e_l ( \frac{1}{2} \delta_{il} k^{- 1} \partial_j ( k ) + \frac{1}{2} \delta_{ij} k^{- 1} \partial_l ( k ) -   \frac{1}{2} \delta_{jl} k^{- 1} \partial_i ( k )   )\\
 &=& \sum_{j,l} e_j \tensora e_l [ ( \Gamma_0 )^i_{jl} + \frac{1}{2} ( \delta_{il} k^{- 1} \partial_j ( k ) +  \delta_{ij} k^{- 1} \partial_l ( k ) -  \delta_{jl} k^{- 1} \partial_i ( k )  )].
		\end{eqnarray*}
		This completes the proof.
		\qed
		
		We complete the section by proving that if $( \E, d )$ is a tame differential calculus such that $\E$ is a free right $\A$-module admitting a central basis, then there indeed exists a bilinear pseudo-Riemannian metric $g_0$ as in Proposition \ref{20thjanprop}. We are going to compute the curvature of the  Levi-Civita connection of such a metric for spectral triples on the noncommutative torus ( Subsection \ref{subnctorus} ) and the quantum Heisenberg manifold ( Subsection  \ref{heisenbergcurvature} ) using the above two results. We will also see that the condition $ d ( e_i ) = 0 $ is satisfied for both these examples. 
		
		\blmma \label{27thmarch20}
				Suppose  $ ( \E, d ) $ is a tame differential calculus such that $ \E $ is a free right $\A$-module with a basis $ \{ e_1, e_2, \cdots e_n \} \in \mathcal{Z} ( \E ). $ Then there exists a unique pseudo-Riemannian bilinear metric $g_0$ on $\E$ such that $ g_0 ( e_i \tensora e_j ) = \delta_{ij}.1. $
		\elmma
		{\bf Proof:} We define 
		\begin{equation} \label{27thmarch202} g_0 : \E \tensora \E \rightarrow \A, ~ g_0 ( ( \sum_i e_i a_i ) \tensora ( \sum_j e_j b_j  ) ) = \sum_i a_i b_i. \end{equation} 
		In particular, $ g_0 ( e_i \tensora e_j ) = \delta_{ij.}$ It is clear that $g_0$ is right $\A$-linear. 
		
		The uniqueness of the map $g_0$ is clear as the facts that  $e_i \in \Ecenter, $  $ g_0 $ is right $\A$-linear and $ g_0 ( e_i \tensora e_j ) = \delta_{ij}.1 $ force $g_0$ to be defined by \eqref{27thmarch202}. The fact that $ g_0 $ is a pseudo-Riemannian metric has been proved in Proposition 2.14 of \cite{article2}. So we only need to check that
		 $g_0$ is bilinear. As remarked above, $ g_0 $ is right $\A$-linear by definition. Let $a$ be an element of $\A.$ Then
		\begin{eqnarray*}
		g_0 (  a ( \sum_i e_i a_i   ) \tensora ( \sum_j e_j b_j  )   ) &=&  g_0 ( ( \sum_i e_i a a_i  ) \tensora ( \sum_j e_j b_j   ) )~ {\rm (} ~ {\rm since} ~ e_i ~ {\rm belong} ~ {\rm to} ~ \Ecenter ~  {\rm )}\\
		&=& \sum_i a a_i b_i  {\rm (} ~ {\rm by} ~ \eqref{27thmarch202} ~ {\rm )}\\
		&=& a ( \sum_i a_i b_i ) \\
		&=& a g_0 ( ( \sum_i e_i a_i   ) \tensora ( \sum_j e_j b_j  )   ) 
		\end{eqnarray*}
		which proves that $g_0$ is bilinear. 
					\qed
		
\section{Computation of the Ricci and scalar curvature} \label{curvaturesection}
		
	Following the groundbreaking work of Connes and Tretkoff in \cite{scalar_3}, computation of  scalar curvature using the asymptotic expansion of the Laplace operator led to several seminal works. We refer to \cite{Connes_moscovici}, \cite{khalkhali}  and references therein. In this section, we take an alternative path. We follow \cite{article3} to compute  the curvature of the Levi Civita connection of a conformally deformed metric on a tame differential calculus. We will apply Proposition \ref{20thjanprop} to compute the Ricci and scalar curvature  for the module of one forms for the canonical spectral triple  on the noncommutative torus. The last subsection will deal with the computation of the curvature for the space of one forms  on the quantum Heisenberg manifold studied in \cite{chak_sinha}.   
	Let us start by defining the notions of Ricci and scalar curvature of a connection on a tame differential calculus $ ( \E, d ). $ For this, we need a few more definitions and some preparatory results. 
	
	Firstly, let us recall ( Definition \ref{tame} ) that the short-exact sequence $ 0 \rightarrow {\rm Ker} ( \wedge ) \stackrel{\iota}{\longrightarrow}  \E \tensora \E \stackrel{\wedge}{\longrightarrow} \twoform $ splits, $\iota$ being the inclusion map. As a result, we have a direct sum decomposition $ \E \tensora \E = {\rm Ker} ( \wedge ) \oplus \mathcal{F}. $ 
	
	Hence, $ \wedge: \F \rightarrow \twoform $ is a right $\A$-linear isomorphism. Moreover, again by definition ( part iv. of Definition \ref{tame}  ), $\Psym$ is an idempotent with range equals to $ {\rm Ker} ( \wedge ) $ and kernel equals to $\F.$ Hence, $ \F  = {\rm Ran} ( 1 - \Psym ).$ Thus, we can view the map $\wedge$ as a right $\A$-linear isomorphism from $ {\rm Ran} ( 1 - \Psym ) $ to $\Omega^2 ( \A ). $ In fact, to avoid confusion about the domain of the map $\wedge,$ we will introduce the following notation.
		
		\bdfn \label{defnQ}
		We will denote the restriction of the map $\wedge$ to $ {\rm Ran} ( 1 - \Psym )$ by the notation $ Q. $
		\edfn
	We make the following observations about the map $Q:$	 
	
		\blmma \label{22ndmarch203}
	For a tame differential calculus $ ( \E, d ) $ and $Q$ as in Definition \ref{defnQ},	we have the following:
		\begin{enumerate}
		\item[i.] The map $ Q: {\rm Ran} ( 1 - \Psym ) \rightarrow \Omega^2 ( \A )$ and its inverse $ Q^{-1}: \twoform \rightarrow {\rm Ran} ( 1 - \Psym ) $ are $\A$-bilinear maps.
		
		\item[ii.] If $ X $ belongs to $\E \tensora \E,$ then 
		   \begin{equation} \label{22ndmarch20} Q (  ( 1 - \Psym  ) ( X )  ) = \wedge ( X ).  \end{equation}
			
		\item[iii.] If $a$ belongs to $\A$ and $f$ belongs to $\E,$ then
			\begin{equation} \label{22ndmarch202} Q^{-1} ( da \wedge f ) = ( 1 - \Psym  ) ( da \tensora f  ). \end{equation}
		\end{enumerate}
		\elmma
		{\bf Proof:} By definition, the map  $Q$ is the restriction of the map $\wedge$ to $ {\rm Ran} ( 1 - \Psym ).$ But by the definition of a differential calculus, the map $\wedge: \E \tensora \E \rightarrow \twoform $ is bilinear  and so $ Q $ is bilinear. Consequently, $ Q^{-1} $ is also  bilinear. 
		
		Next, if $X$ belongs to $\E \tensora \E,$  the equation \eqref{22ndmarch20} follows from the following computation:
	$$ Q (  ( 1 - \Psym ) ( X ) ) = \wedge ( ( 1 - \Psym  ) ( X )  ) = \wedge ( X ) - \wedge ( \Psym ( X )  ) = \wedge ( X ) - 0  $$
	as $ \Psym ( X )  $ belongs to $ {\rm Ran} ( \Psym ) $ which is equal to $ {\rm Ker} ( \wedge ) $ by part iv. of Definition \ref{tame}.
		
		Now we use \eqref{22ndmarch20} to prove \eqref{22ndmarch202}. Since $ (  1 - \Psym ) ( X ) \in {\rm Ran} ( 1 - \Psym )$ and $Q: {\rm Ran} ( 1 - \Psym ) \rightarrow \twoform $ is a right $\A$-linear isomorphism, \eqref{22ndmarch20} allows us to conclude that for all $X \in \E \tensora \E,$ 
	\begin{equation} \label{22ndmarch205} Q^{-1} ( \wedge ( X ) ) = ( 1 - \Psym ) ( X ). \end{equation}
			In particular, if $a$ belongs to $\A$ and $f$ belongs to $\E,$ then
			$$ Q^{-1} ( da \wedge f ) = Q^{-1} ( \wedge ( da \tensora f )  ) = ( 1 - \Psym   ) ( da \tensora f  ).$$
		This completes the proof of the lemma.
		\qed
		
		We will need another lemma to define the curvature operator.
		
		\blmma \label{10thjuly2018}
		Let	$ ( \E, d ) $ be a tame differential calculus and $\nabla$ a connection on $\E.$ Then the map
		$$ \E \tensorc \E \rightarrow \E \tensora \E \tensora \E ~ {\rm defined} ~ {\rm by} ~ e \tensorc f \mapsto ( 1 - \Psym   )_{23} ( \nabla e \tensora f ) + e \tensora Q^{ - 1 } ( d f ) $$
			descends  to a  map from  $ \E \tensora \E $ to $ \E \tensora \E \tensora \E. $ We will denote this map by the symbol $H.$
			
Moreover, if we define 
$$ R ( \nabla ):= H \circ \nabla : \E \rightarrow \E \tensora \E \tensora \E, $$
				then $ R ( \nabla ) $ is a right $ \A $-linear map.			
\elmma	
			{\bf Proof:}	
			Let  $e, f \in \E $ and $ a \in \A. $ Then we get 
					\begin{eqnarray*}
					 H ( e \tensorc a f ) &=& ( 1 - \Psym )_{23} ( \nabla ( e ) \tensora a f ) + e \tensora Q^{-1} ( da \wedge f + a. d f   )\\
					&=& ( 1 - \Psym )_{23} ( \nabla ( e ) \tensora a f ) +  e \tensora ( 1- \Psym ) ( da \tensora f ) + e \tensora Q^{-1} (   a d f   ) \\
					&& {\rm (} ~ {\rm by} ~ \eqref{22ndmarch202} ~ {\rm )}\\
					&=& ( 1 - \Psym )_{23} ( \nabla ( e ) a \tensora  f ) + e \tensora ( 1- \Psym ) ( da \tensora f )  +  e a \tensora Q^{-1} ( d f ) \\
					&& {\rm (} ~ {\rm since}  ~ Q ~ {\rm is} ~ {\rm left} ~ \A-{\rm linear} ~ {\rm by} ~ Lemma ~ \ref{22ndmarch203} ~ {\rm )}\\
					&=& ( 1 - \Psym )_{23} ( ( \nabla ( e ) a + e \tensora da ) \tensora  f ) + e a \tensora Q^{-1} ( d f )\\
					&& {\rm(} ~ {\rm as} ~ 1 - \Psym ~ {\rm is} ~ {\rm an} ~ {\rm idempotent} ~ {\rm )}\\
					&=& ( 1 - \Psym )_{23} (  \nabla ( e a )  \tensora f ) + e a \tensora Q^{-1} ( d f )\\
					&=& H ( e a \tensorc  f ),
					\end{eqnarray*}
					which proves that $H$ descends  to a  map from  $ \E \tensora \E $ to $ \E \tensora \E \tensora \E. $
				
				Now we prove that the map $R ( \nabla ): \E \rightarrow \E \tensora \E \tensora \E $ is a right $\A$-linear map. Let $ e \in \E $ and $ a \in \A. $ We will use the Sweedler-type notation $ \nabla ( e ) = e_{(1)} \tensora e_{(2)}. $ It follows that 
\begin{equation} \label{22ndmarch204} R ( \nabla ) ( e ) = ( 1 - \Psym  )_{23} (  \nabla ( e_{(1)} ) \tensora e_{(2)}   ) + e_{(1)} \tensora Q^{-1} (  d ( e_{(2)}  ) ). \end{equation}
				Therefore, using  $ \nabla ( e a ) = \nabla ( e ) a + e \tensora da, $   it is easy to see that
				\begin{eqnarray*}
				  R ( \nabla ) ( e a ) &=& ( 1 - \Psym  )_{23} (  \nabla ( e_{(1)} ) \tensora e_{(2)} a   ) + e_{(1)} \tensora Q^{-1} (  d ( e_{(2)} a  ) ) + ( 1 - \Psym  )_{23} (  \nabla ( e ) \tensora da   )\\
					&=& ( 1 - \Psym  )_{23} (  \nabla ( e_{(1)} ) \tensora e_{(2)}    ) a + e_{(1)} \tensora [ Q^{-1} (  d ( e_{(2)}  ) a - e_{(2)} \wedge da  ) + ( 1 - \Psym  ) (  e_{(2)} \tensora da  )    ]\\
					&& {\rm (} ~ {\rm since} ~ \Psym ~ {\rm and} ~ Q^{-1} ~ {\rm are} ~ {\rm right} ~ \A-{\rm linear},~ {\rm see} ~ {\rm Lemma} ~ \ref{22ndmarch203} ~ {\rm )}\\
					&=& ( 1 - \Psym  )_{23} (  \nabla ( e_{(1)} ) \tensora e_{(2)}   ) a + e_{(1)} \tensora Q^{-1} (  d ( e_{(2)}  ) ) a\\
					&& {\rm(} ~ {\rm by} ~ \eqref{22ndmarch205} ~ {\rm )}\\
					&=& R ( \nabla ) ( e ) a
				\end{eqnarray*}
				by \eqref{22ndmarch204}, proving that $ R ( \nabla ) $ is right $\A$-linear. This finishes the proof of the lemma.
				\qed
				
		Now we are prepared to define  the ``curvature operator" following \cite{article3} and \cite{chak_sinha}. We observe that since $ R ( \nabla )  $ belongs to $ \Hom_{\A} ( \E, \E \tensora \E \tensora \E ) $ by Lemma \ref{10thjuly2018}, we can apply the map $ \zeta^{-1}_{\E, \E \tensora \E \tensora \E} $ ( see Definition \ref{xi} ) to $ R ( \nabla ) $ and the image lies in $ ( \E \tensora \E \tensora \E ) \tensora \E^*. $
		
		\bdfn \label{23rdmarch20}
		If $( \E, d )$ is a tame differential calculus and $\nabla$ is a torsionless connection on $\E,$ the curvature operator $ \Theta $ of the connection $\nabla$ is defined to be the image of the element $ R ( \nabla ) $ under the following maps:
		\[ 
 \begin{tikzcd}
  \Hom_{\A} ( \E, \E \tensora \E \tensora \E ) \arrow[r, "\zeta^{-1}_{\E, \E \tensora  \E \tensora  \E}"] & [7 ex] ( \E \tensora  \E \tensora \E ) \tensora \E^* \arrow[r, "\sigma_{23}"] &  \E \tensora \E  \tensora \E \tensora \E^*.
 \end{tikzcd}
\]
Here, $ \sigma_{23}:\E \tensora \E \tensora \E \rightarrow \E \tensora \E \tensora \E $ is the map $  {\rm id}_\E \tensora \sigma.$
	\edfn	
		
		Now we proceed towards the definitions of the Ricci curvature and scalar curvature. We will need a lemma whose proof is elementary:
		
		\blmma
		If $( \E, d ) $ is a tame differential calculus, $ u^\E: \Ecenter \otimes_{\Acenter} \A \rightarrow \E $ be the multiplication map defined in part iii. of Definition \ref{tame}, $  v^\E: \A \otimes_{\mathcal{Z} ( \A )} \mathcal{Z} ( \E ) \rightarrow \E  $ be defined by
	 $$ v^\E \left( \sum_i a_i \otimes_{\mathcal{Z} ( \A )} \omega_i \right)  = \sum_i a_i \omega_i ~ {\rm and}$$
	$$ {\rm flip}: \Ecenter \otimes_{\Acenter} \E^* \rightarrow \E^* \otimes_{\Acenter} \Ecenter ~ {\rm defined} ~ {\rm by} ~  {\rm flip} ~ ( e^\prime \otimes_{\Acenter} \phi ) = \phi \otimes_{\Acenter} e^\prime $$
	for all $ e^\prime $ in $\Ecenter$ and $\phi$ in $\E^*,$  then the map $ \rho: \E \tensora \E^* \rightarrow \E^* \tensora \E $ defined 
			as the composition:
			\[ 
 \begin{tikzcd}
  \E \tensora \E^* \arrow[r, "( u^{\E} )^{-1} \tensora {\rm id}_{\E^*}"] & [5 ex] \Ecenter \otimes_{\Acenter} \E^* \arrow[r, "{\rm flip}"] & [5 ex]  \E^* \otimes_{\Acenter} \Ecenter = \E^* \tensora ( \A \otimes_{\Acenter} \Ecenter  ) \arrow[r, "{\rm id}_{\E^*} \tensora v^{\E}"] & [5 ex] \E^* \tensora \E.
 \end{tikzcd}
\]
	is actually well-defined.	
		\elmma
		{\bf Proof:} We only need to check that each of the maps $ (  u^{\E}  )^{-1} \tensora {\rm id}_{\E^*}: \E \tensora \E^* \rightarrow \Ecenter \otimes_{\Acenter} \E^*, ~ {\rm id}_{\E^*} \tensora v^{\E}: \E^* \otimes_{\Acenter} \Ecenter \rightarrow \E^* \tensora \E $ and $ {\rm flip}:  \Ecenter \otimes_{\Acenter} \E^* \rightarrow \E^* \otimes_{\Acenter} \Ecenter $ is well-defined.

		By  Remark \ref{17thdec2019remark}, we know  that the map $u^{\E}  $ defined as
		$$ u^{\E} ( \sum_i \omega_i \otimes_{\Acenter} a_i ) = \sum_i \omega_i a_i $$
		is a right $\A$-linear isomorphism. Thus, the map
		$ (  u^{\E}  )^{-1} \tensora {\rm id}_{\E^*}: \E \tensora \E^* \rightarrow \Ecenter \otimes_{\Acenter} \E^* $ is well-defined.
		
		Moreover, an inspection of the proof of Proposition 2.4 of \cite{article1} shows that $ v^{\E}: \A \otimes_{\Acenter} \Ecenter \rightarrow \E $ is a left $\A$-linear, right $\Acenter$-linear and invertible.  
	Thus, the map $ {\rm id}_{\E^*} \tensora v^{\E}: \E^* \otimes_{\Acenter} \Ecenter  = \E^* \tensora ( \A \otimes_{\Acenter} \Ecenter  ) \rightarrow \E^* \tensora \E $ is well-defined.

		Finally, it can be easily checked that $ {\rm flip} $ is well-defined ( and a right $ \Acenter $-linear isomorphism ).
			\qed

	Now we are prepared to define the Ricci curvature and the scalar curvature.

			\bdfn
			For a tame differential calculus $ ( \E, d ) $ and a torsionless connection	 $\nabla$ on $\E,$ the Ricci curvature  $ {\rm Ric} $  is defined as the element in $ \E \tensora \E $ given by
				\be \label{23rdjan1} {\rm Ric}:= ( {\rm id}_{\E \tensora \E} \tensora {\rm ev} \circ \rho ) ( \Theta ),\ee
				where $ {\rm ev}: \E^* \tensora \E \rightarrow \A $ is the $\cla $-bilinear map sending $ e^* \tensora f $ to $ e^* ( f ) $ for all $ e^* \in \E^* $ and $ f \in \E $ and $\Theta$ is the curvature operator defined in Definition \ref{23rdmarch20}.
				
				Finally, the scalar curvature {\rm Scal} is defined as:
				\be \label{23rdjan2} {\rm Scal}: = {\rm ev} ( V_g \tensora {\rm id}_\E  ) ( {\rm Ric} ) ~ \in ~ \cla. \ee
				\edfn	
					
					\brmrk It is easy to see that in the classical case, i.e, when $ \E = \oneform $ and $\A = C^\infty ( M ), $ the above definitions of Ricci and scalar curvature do coincide with the usual notions.
					\ermrk

					\bppsn \label{20thjanprop0}
					If  $ ( \E, d ) $ is a tame differential calculus such that $\E$ is a free right $\A$-module with a basis $ \{ e_1, e_2, \cdots e_n \} \in \mathcal{Z} ( \E ) $ such that   $ d ( e_i ) = 0 $ for all $ i = 1,2, \cdots n.$ Then the curvature operator, Ricci tensor and the scalar curvature of a torsion-less connection $\nabla$ are given by:
				$$ R ( \nabla ) ( e_i ) = \sum_{j,k,l} e_j \tensora e_k \tensora e_l r^i_{jkl} ~ {\rm where} ~ r^i_{jkl} = \frac{1}{2} \sum_p [ ( \Gamma^p_{jk} \Gamma^i_{pl} - \Gamma^p_{jl} \Gamma^i_{pk} ) - \partial_l ( \Gamma^i_{jk} ) + \partial_k ( \Gamma^i_{jl} )  ].$$
					The Ricci tensor $ {\rm Ric} $ is given by $ {\rm Ric} = \sum_{j,l} e_j \tensora e_l {\rm Ric} ( e_j, e_l ) $ where
 $$ {\rm Ric} ( e_j, e_l ) = \frac{1}{2} \sum_i [ \sum_p ( \Gamma^p_{ji} \Gamma^i_{pl} - \Gamma^p_{jl} \Gamma^i_{pi} ) - \partial_l ( \Gamma^i_{ji} ) + \partial_i ( \Gamma^i_{jl} ) ].  $$
  The Scalar curvature is given by $ {\rm Scal} = \sum_{j,l} g ( e_j \tensora e_l ) {\rm Ric} ( e_j, e_l ).$
	\eppsn
	
	{\bf Proof:} We use the definition of $H$ in Lemma \ref{10thjuly2018} to compute
	  \begin{eqnarray}
		R ( \nabla ) ( e_i ) &=& H \circ \nabla ( e_i ) \nonumber \\
		                     &=& H ( \sum_{j,k} e_j \tensora e_k \Gamma^i_{jk}   ) \nonumber \\
												 &=& \sum_{j,k} [  ( 1 - \Psym  )_{23} (  \nabla ( e_j ) \tensora e_k \Gamma^i_{jk}  ) + e_j \tensora Q^{-1} (  d ( e_k \Gamma^i_{jk}  )  )  ] \nonumber \\
												&=& \sum_{j,k} [  \frac{1 - \sigma_{23}}{2} (  \sum_{m,n} e_m \tensora e_n \Gamma^j_{mn} \tensora e_k \Gamma^i_{jk}   )+ e_j \tensora Q^{-1} (  d ( e_k ) \Gamma^i_{jk} - e_k \wedge d ( \Gamma^i_{jk} )  )   ] \nonumber \\
												&=& \sum_{j,k,m,n} \frac{1 - \sigma_{23}}{2} ( e_m \tensora e_n \tensora e_k \Gamma^j_{mn} \Gamma^i_{jk}  ) + \sum_{j,k} e_j \tensora Q^{-1} ( d ( e_k )  ) \Gamma^i_{jk} \nonumber \\ 
												&-& \sum_{j,k} e_j \tensora Q^{-1} ( e_k \wedge d (  \Gamma^i_{jk}  )   ) \label{24thmarch20}
		\end{eqnarray}
since $e_i$ belongs to $\Ecenter.$
Now,
\begin{eqnarray}
&& \sum_{n,k} \frac{1 - \sigma_{23}}{2} ( e_m \tensora e_n \tensora e_k \Gamma^j_{mn} \Gamma^i_{jk}    ) \nonumber \\
&=& \frac{1}{2} \sum_{n,k} e_m \tensora e_n \tensora e_k ( \Gamma^j_{mn} \Gamma^i_{jk} - \Gamma^j_{mk} \Gamma^i_{jn}   ) \label{24thmarch202}
\end{eqnarray} 	
by applying \eqref{10thjuly20182} as $e_n, e_k$ belong to $\Ecenter.$
Next, as $\nabla$ is a torsionless connection, we have
\begin{eqnarray*}
Q^{-1} ( d ( e_k ) ) &=& Q^{-1} ( - \wedge \circ \nabla ( e_k ) ) \\
                     &=& - ( 1 - \Psym ) \nabla ( e_k ) ~ {\rm (} ~ {\rm by} ~ \eqref{22ndmarch205} ~ {\rm )} \\
										 &=& - \frac{1 - \sigma}{2} ( \sum_{r,s} e_r \tensora e_s \Gamma^k_{rs}  ) \\
										 &=& - \frac{1}{2} \sum_{r,s} e_r \tensora e_s ( \Gamma^k_{rs} - \Gamma^k_{sr}  )
\end{eqnarray*}
by applying \eqref{10thjuly20182}. However, $ \Gamma^k_{rs} = \Gamma^k_{sr} $ by \eqref{20thjan2} and hence
\begin{equation} \label{24thmarch203} Q^{-1} ( d ( e_k ) ) = 0. \end{equation}
Finally, by applying \eqref{22ndmarch205}, we get
\begin{eqnarray}
&& \sum_k Q^{-1} ( e_k \wedge d ( \Gamma^i_{jk}  ) ) \nonumber \\
&=& \sum_k ( 1 - \Psym ) ( e_k \tensora d ( \Gamma^i_{jk} ) \nonumber \\
&=& \frac{1}{2} \sum_{k,n} ( 1 - \sigma ) ( e_k \tensora e_n \partial_n (  \Gamma^i_{jk} )   ) ~ {\rm (} ~ {\rm by} ~ \eqref{19thjan3} ~ {\rm )} \nonumber \\
&=& \frac{1}{2} \sum_{k,n} e_k \tensora e_n ( \partial_n (  \Gamma^i_{jk} ) - \partial_k ( \Gamma^i_{jn} )  ) ~ {\rm (} ~ {\rm by} ~ \eqref{10thjuly20182} ~ {\rm )}. \label{24thmarch204} 
\end{eqnarray} 

Plugging \eqref{24thmarch202}, \eqref{24thmarch203} and \eqref{24thmarch204} in \eqref{24thmarch20}, we obtain
\begin{eqnarray*}
&& R ( \nabla ) ( e_i )\\
&=& \frac{1}{2} \sum_{j,k,m,n} e_m \tensora e_n \tensora e_k ( \Gamma^j_{mn} \Gamma^i_{jk} - \Gamma^j_{mk} \Gamma^i_{jn} ) - \frac{1}{2} \sum_{j,k,n} e_j \tensora e_k \tensora e_n ( \partial_n ( \Gamma^i_{jk} ) - \partial_k (  \Gamma^i_{jn}  )  )\\
&=& \sum_{j,k,l} e_j \tensora e_k \tensora e_l [  \frac{1}{2} \sum_p ( \Gamma^p_{jk} \Gamma^i_{pl} - \Gamma^p_{jl} \Gamma^i_{pk} )  ] + \sum_{j,k,l} e_j \tensora e_k \tensora e_l . \frac{1}{2} [ - \partial_l (  \Gamma^i_{jk}  ) + \partial_k ( \Gamma^i_{jl}  ) ]\\
&=& \sum_{j,k,l} e_j \tensora e_k \tensora e_l [ \frac{1}{2} \left( \sum_p ( \Gamma^p_{jk} \Gamma^i_{pl} - \Gamma^p_{jl} \Gamma^i_{pk}  ) - \partial_l (  \Gamma^i_{jk}  ) + \partial_k (  \Gamma^i_{jl}  ) \right) ].
\end{eqnarray*}
This completes the proof of the result.
\qed
	
	\brmrk
	The only place where we have used the condition $ d ( e_i ) = 0 $ is the equality $ \Gamma^k_{rs} = \Gamma^k_{sr} $ which proves that $ Q^{-1} ( d ( e_k ) ) = 0 $ ( \eqref{24thmarch203} ). Thus, in the absence of this condition, we would get some additional terms.   				
	\ermrk

		\subsection{Computation of curvature for the conformally deformed metric on the noncommutative torus} \label{subnctorus}
		
We recall  that the noncommutative 2-torus $ C (  \IT^2_\theta ) $ is the universal $ C^* $ algebra generated by two unitaries $ U $ and $ V $ satisfying $ U V = e^{2 \pi i \theta} V U $ where $ \theta $ is a number in $ [ 0, 1 ].$  The $\ast$- subalgebra $ \cla ( \IT^2_\theta ) $ of $  C (  \IT^2_\theta ) $ generated by $ U $ and $ V $ will be denoted by $ \cla. $

		We  have the following concrete description of a spectral geometry on $ \cla $:  (  see \cite{connes_landi}   ):
		
 there are two derivations $ d_1 $ and $ d_2 $ on $ \cla  $ obtained by extending linearly the rule: 
      $$ d_1 ( U ) = U, ~ d_1 ( V ) = 0,~ d_2 ( U ) = 0, ~ d_2 ( V ) = V. $$  
       There is a  faithful trace on $ \cla  $ defined as follows:
 $$ \tau ( \sum_{m,n} a_{mn} U^m V^n ) = a_{00}, ~ {\rm where} ~ {\rm the} ~ {\rm sum} ~ {\rm runs} ~ {\rm over} ~ {\rm a} ~  {\rm finite} ~ {\rm subset} ~ {\rm of} ~ \IZ \times \IZ. $$
 Let $ \clh = L^2 ( C ( \IT^2_\theta ), \tau ) \oplus L^2 ( C (  \IT^2_\theta  ), \tau ) $ where $ L^2 ( C (   \IT^2_\theta  ),  \tau ) $ denotes the GNS Hilbert space of $ \cla  $ with respect to the state $ \tau.$ We note that $ \cla  $ is embedded as a subalgebra of $ \clb ( \clh ) $ by $ a \mapsto \left(  \begin {array} {cccc}
     a   &  0  \\ 0 & a \end {array} \right ) .$
     The Dirac operator on $ \clh $ is defined by 
		$$ D = \left(  \begin {array} {cccc}
     0   & d_1 + \sqrt{- 1} d_2  \\ d_1 - \sqrt{- 1} d_2 & 0 \end {array} \right ).$$
     Then, $( \cla, \clh, D )  $ is a spectral triple of compact type.  In particular, for $ \theta = 0, $ this coincides with the classical spectral triple on $ \IT^2. $
			
			Let $ \gamma_1 = \left(  \begin {array} {cccc}
     0   & 1  \\ 1 & 0 \end {array} \right ) $ and $ \gamma_2 = \left(  \begin {array} {cccc}
     0   & \sqrt{- 1}  \\ - \sqrt{- 1} & 0 \end {array} \right ) .$
		The de-Rham differential $ d:= d_D  :  \cla \rightarrow \E:= \oneform  $ is defined by
		 $$ d ( a ) = \sqrt{- 1} [ D, a ]. $$
		
		We have the following result:
		
		\bppsn \label{22ndjan}
		The differential calculus $ ( \E, d ) $ is tame. In fact, 
		the bimodule $ \E $ of one-forms is freely generated as a right $\A$-module by the central elements
		$$ e_1 = 1 \tensorc \gamma_1, ~ e_2 = 1 \tensorc \gamma_2.$$
		The space of two forms is a rank one free module generated by $e_1 \wedge e_2.$
		Moreover, we have  
		\be \label{7thmay2018} d ( e_1 ) = d ( e_2 ) = 0.\ee   
\eppsn
		{\bf Proof:} Consider the usual spectral triple on the $2$-torus. Then the group $ \mathbb{T}^2 $ acts freely and isometrically on $ \mathbb{T}^2.$  The spectral triple on $\A$ defined above is the isospectral deformation ( \cite{connes_landi} ) of the classical spectral triple and hence we can apply Theorem 7.1 of \cite{article1} to conclude that $ ( \E, d ) $ is tame.
		The structure of the space of one forms and two forms is well known and hence we omit the proof. It is also clear that $e_1$ and $e_2$ are  elements of $\mathcal{Z} ( \E ).$ 
		
				Thus, we are left with proving  \eqref{7thmay2018}. 
Since $\gamma_1. \gamma_2 = - \gamma_2. \gamma_1 $ and $ \gamma^2 = 1 $ we have 
		\begin{equation} \label{26thmarch20} e_1 \wedge e_1 = 0 ~ {\rm and} ~ e_1 \wedge e_2 = - e_2 \wedge e_1. \end{equation}  
		Now,	it is easy to see that 
		$$ d U = \sqrt{- 1} e_1 U, ~ d V = \sqrt{- 1} e_2 V. $$ 
Therefore, by  Leibniz rule, we have 
$$ 0 = d^2 U =  d ( e_1 ). U - e_1 \wedge d U  = d ( e_1 ) U -  \sqrt{- 1}  e_1 \wedge e_1. U.  $$
By \eqref{26thmarch20}, we obtain 
  $ d ( e_1 ) U = 0 $ and hence $ d ( e_1 ) = 0 $  since $U$ is invertible. Similarly, $ d ( e_2 ) = 0. $ \qed

\brmrk \label{24thjan1}
Proposition \ref{22ndjan} allows us to apply Proposition \ref{20thmarch203} to the differential calculus $(\E, d ).$ From the equalities $ d U = \sqrt{- 1} e_1 U, ~ d V =  \sqrt{- 1} e_2 V,$ it follows that the derivations $ \partial_1 $ and $ \partial_2 $ as in \eqref{19thjan3} are given by the following formulas:
  $$ \partial_1 ( U ) = \sqrt{- 1} U, ~ \partial_1 ( V ) = 0, ~ \partial_2 ( U )  = 0, \partial_2 ( V ) = \sqrt{- 1} V. $$
		From these formulas, it can be easily checked that  $ \partial_1 $ and $ \partial_2 $ commute.
 \ermrk

\bthm \label{1stoct204}
Consider the differential calculus $ ( \E, d ) $ on the noncommutative $2$-torus $\A$ as above. Consider the pseudo-Riemannian bilinear metric $g_0$ of Lemma \ref{27thmarch20}. Let $ k $ be an invertible element of $\A.$ Then the Ricci and the scalar curvatures of the Levi-Civita connection for the triplet $ ( \E, d, k. g_0 ) $  are as follows: 
$$ {\rm Ric} ( e_1, e_1 ) = {\rm Ric} ( e_2, e_2 ) = - \frac{1}{2} (  k^{- 1} ( \partial^2_1 + \partial^2_2 ) ( k ) + \partial_1 ( k^{ - 1 }  ) \partial_1 ( k ) + \partial_2 ( k^{ - 1 } ) \partial_2 ( k ) ).$$
$$ {\rm Ric} ( e_1, e_2 ) = -  {\rm Ric} ( e_2, e_1 ) = \frac{1}{2} ( \partial_1 ( k^{ - 1 } ) \partial_2 ( k ) - \partial_2 ( k^{ - 1 }  ) \partial_1 ( k ) ). $$
$$ {\rm Scal} = - ( \partial^2_1 + \partial^2_2  ) ( k ) - k ( \partial_2 ( k^{ - 1 } ) \partial_2 ( k ) - k \partial_1 ( k^{ - 1 }  ) \partial_1 ( k ) ).$$
\ethm
{\bf Proof:} By virtue of  Proposition \ref{22ndjan}, the hypotheses of Lemma \ref{27thmarch20} and Proposition \ref{20thjanprop0} hold and so we can apply them to this differential calculus. In particular, 
\be \label{20thjan} g_0  ( e_i \tensora e_j ) = \delta_{ij} 1_{\cla} \ee
and Theorem \ref{19thjan1} ensures the existence and uniqueness of the Levi-Civita connection for the triplet $ ( \E, d, k. g_0 ) $ with Christoffel symbols as in \eqref{22ndjan2}.

We will use Proposition \ref{20thjanprop} to compute the Christoffel symbols of the Levi-Civita connection for the triplet $( \E, d, g ).$ Let $ \nabla_{g_0} $ be the Levi-Civita connection for the triplet $( \E, d, g_0 )$ as in Proposition \ref{20thjanprop}. We claim that 
\begin{equation} \label{30thmarch20} \nabla_{g_0} ( e_i ) = 0 ~ {\rm for} ~ {\rm all} ~ i. \end{equation}
Indeed, we define a connection $\nabla_1$ on $\E$ by the formula
\begin{equation} \label{30thmarch202} \nabla_1 ( \sum_i e_i a_i ) = \sum_i e_i \tensora d a_i. \end{equation}
In particular, $ \nabla_1 ( e_i ) = 0 $ for all $i.$

We prove that $ \nabla_1 $ is a torsionless and compatible with $g_0$ so that the uniqueness of the Levi-Civita connection for a {\bf bilinear} pseudo-Riemannian metric ( Theorem \ref{21stfeb2020} ) will imply that $ \nabla_1 = \nabla_{g_0}. $ We have
\begin{eqnarray*}
 \wedge \circ \nabla_1 ( \sum_i e_i a_i ) &=& \sum_i \wedge ( e_i \tensora da ) ~ {\rm (} ~ {\rm by} ~ \eqref{30thmarch202} ~ {\rm )}\\
                                          &=& \sum_i e_i \wedge  d a_i\\
																					&=& - d ( \sum_i e_i a_i )
																					\end{eqnarray*} 
as $ d ( e_i ) = 0 $ (  Proposition \ref{22ndjan} ) proving that $\nabla_1$ is torsionless. 

Next, by Definition \ref{19thfeb202021}, $\nabla_1$ is compatible with $g_0$ if for all $e, f$ in $\E,$
$$ \Pi_{g_0} ( \nabla_1 ) ( e \tensora f ) = d ( g ( e \tensora f )  ).  $$
By \eqref{20thfeb2020}, for all $a $ in $\A,$
we have
$$ \Pi_{g_0} ( \nabla_1 ) ( \sum_{i,j} e_i \tensora e_j a  ) = \sum_{i,j} \Pi^0_{g_0} ( \nabla_1 ) ( e_i \tensora e_j ) a + \sum_{i,j} g_0 ( e_i \tensora e_j ) da,  $$
where, from Proposition \ref{30thmarch203},
$$ \Pi^0_{g_0} ( \nabla_1 ) ( e_i \tensora e_j ) = ( g_0 \tensora {\rm id}  ) \sigma_{23} ( \nabla_1 ( e_i ) \tensora e_j + \nabla_1 ( e_j ) \tensora e_i ) = 0 $$
as $ \nabla_1 ( e_i ) = 0 $ for all $i.$ Therefore, 
\begin{eqnarray*}
\Pi_{g_0} ( \nabla_1 ) ( \sum_{i,j} e_i \tensora e_j a  ) &=& 0 + \sum_{i,j} g_0 ( e_i \tensora e_j  ) da\\
                                                          &=& \sum_{i,j} [ d ( g_0 ( e_i \tensora e_j )   ) a + g_0 ( e_i \tensora e_j ) da  ] ~ {\rm (} ~ {\rm as} ~ g_0 ( e_i \tensora e_j  ) = \delta_{ij}.1 ~ {\rm by} ~ \eqref{20thjan} ~   {\rm )}\\
																													&=& \sum_{i,j} d ( g_0 ( e_i \tensora e_j a )  ) ~ {\rm (} ~ {\rm as} ~ g_0 ~ {\rm is} ~ {\rm right} ~ \A-{\rm linear} ~ {\rm and} ~ {\rm we} ~ {\rm have} ~ {\rm appplied} ~ {\rm Leibniz} ~ {\rm rule} ~ {\rm )}\\
																													&=& d ( g_0 ( \sum_{i,j} e_i \tensora e_j a  )  ).
\end{eqnarray*}
Hence, $\nabla_1$ is compatible with $g_0.$ By the discussion made above, $\nabla_1 = \nabla_{g_0} $ and so \eqref{30thmarch20} holds. By applying \eqref{22ndjan2}, we obtain 
  $$ \Gamma^i_{jl} = \frac{1}{2} ( \delta_{il} k^{- 1} \partial_j ( k ) +  \delta_{ij} k^{- 1} \partial_l ( k ) -  \delta_{jl} k^{- 1} \partial_i ( k )  ),$$
 $ \partial_i $ being the derivations as in \eqref{19thjan3}.																					
 Thus, we have:
$$ \Gamma^1_{11} = \frac{1}{2} k^{- 1} \partial_1 ( k ), ~ \Gamma^1_{22} = - \frac{1}{2} k^{ - 1 } \partial_1 ( k ), ~ \Gamma^1_{12} = \Gamma^1_{21} = \frac{1}{2} k^{ - 1 } \partial_2 ( k ), $$
$$ \Gamma^2_{11} = - \frac{1}{2} k^{ - 1 } \partial_2 ( k ), ~ \Gamma^2_{22} = \frac{1}{2} k^{ - 1 } \partial_2 ( k ), ~ \Gamma^2_{12} = \Gamma^2_{21} = \frac{1}{2} k^{ - 1} \partial_1 ( k ). $$

Using these formulas for Christoffel symbols,  we can compute ( using  Proposition \ref{20thjanprop0} ),
\begin{eqnarray*} {\rm Ric} ( e_1, e_1 ) &=& \sum^2_{i,p = 1} ( \Gamma^p_{1i} \Gamma^i_{p1} - \Gamma^p_{11} \Gamma^i_{pi} ) - \sum^2_{i = 1} ( \partial_1 ( \Gamma^i_{1i} ) - \partial_i ( \Gamma^i_{11} ) )\\
 &=& \Gamma^1_{12} \Gamma^2_{11} - \Gamma^1_{11} \Gamma^2_{12} + \Gamma^2_{12} \Gamma^2_{21} - \Gamma^2_{11} \Gamma^2_{22} - \frac{1}{2} ( \partial_1 ( k^{ - 1 } \partial_1 ( k ) ) + \partial_2 ( k^{ - 1 } \partial_2 ( k ) )  )\\
&=& 0 - \frac{1}{2} \partial_1 ( k^{ - 1 } ) \partial_1 ( k ) - \frac{1}{2} k^{ - 1 } \partial^2_1 ( k ) -  \frac{1}{2} \partial_2 ( k^{ - 1 } ) \partial_2 ( k ) -  \frac{1}{2} k^{ - 1 } \partial^2_2 ( k )\\
&=& - \frac{1}{2} (  k^{- 1} ( \partial^2_1 + \partial^2_2 ) ( k ) + \partial_1 ( k^{ - 1 }  ) \partial_1 ( k ) + \partial_2 ( k^{ - 1 } ) \partial_2 ( k ) ).
\end{eqnarray*} 
The computations for $ {\rm Ric} ( e_2, e_2 ), ~ {\rm Ric} ( e_1, e_2 ) $ and $ {\rm Ric} ( e_2, e_1 ) $ are similar and hence omitted. The only extra ingredient   in the computation of $ {\rm Ric} ( e_1, e_2 ) $ and $ {\rm Ric} ( e_2, e_1 ) $ is that the derivations $ \partial_1 $ and $ \partial_2 $ commute as was remarked in Remark \ref{24thjan1}.

Finally, using the formula of the scalar curvature in   Proposition \ref{20thjanprop0} and the equation \eqref{20thjan}, we get that
\begin{eqnarray*}
 {\rm Scal} &=& \sum_{j,l} k g_0 ( e_j \tensora e_l ) {\rm Ric} ( e_j, e_l )\\
            &=& \sum_j k {\rm Ric} ( e_j, e_j )\\
						&=&  - k ( k^{ - 1 } ( \partial^2_1 + \partial^2_2  ) ( k ) + \partial_2 ( k^{ - 1 } ) \partial_2 ( k ) + \partial_1 ( k^{ - 1 }  ) \partial_1 ( k ) )\\
						&=& - ( \partial^2_1 + \partial^2_2  ) ( k ) - k ( \partial_2 ( k^{ - 1 } ) \partial_2 ( k ) - k \partial_1 ( k^{ - 1 }  ) \partial_1 ( k ) ).
						\end{eqnarray*} \qed 
						
	In subsection 3.1 of \cite{article1}, a canonical candidate $g^\prime$ for a pseudo-Riemannian bilinear metric on a tame spectral triple was proposed. It can be easily checked that for the spectral triple on the noncommutative torus under consideration, $g^\prime$ is indeed a pseudo-Riemannian bilinear metric. The proof follows along the lines of Proposition 6.4 of \cite{article1}. It can be easily seen that $ g^\prime ( e_i \tensora e_j ) = g_0 ( e_i \tensora e_j ) = \delta_{ij} $ and so by the uniqueness of Lemma \ref{27thmarch20}, $ g^\prime = g_0. $

						\subsection{Computation of the curvature for the example of the quantum Heisenberg manifold} \label{heisenbergcurvature}
						
						In this subsection, we compute the curvature of the Levi Civita  connection for a certain metric on the space of one-forms of the quantum Heisenberg manifold. The $C^*$-algebra of the quantum Heisenberg manifold was defined and studied in \cite{rieffel_heisenberg}.   The differential calculus which we will consider comes from a spectral triple constructed in \cite{chak_sinha}. For the precise definition of the algebra of the quantum Heisenberg manifold and the spectral triple on it, we refer to \cite{chak_sinha}.  The authors of \cite{chak_sinha} proved that there exists a pseudo-Riemannian metric on the space of one forms for which there is no torsion-less connection which is also metric compatible in the sense of \cite{frolich}. However, using our definition of metric compatibility of a connection, it has been proved in \cite{article1} ( Theorem 6.4 ) that there exists a unique Levi-Civita connection for any pseudo-Riemannian bilinear metric.

				In the rest of this subsection, we will be using the notations and results  of Section 6 of \cite{article1} as well as \cite{chak_sinha}. In particular, we have the following:
				\bppsn ( \cite{chak_sinha}, \cite{article1} ) \label{30thmarch20n}
		Let $ ( \E, d ) $ denote the differential calculus on the quantum Heisenberg manifold $\A$ as in \cite{chak_sinha}. Then the bimodule of one forms $\E$	is a free right $\A$-module of rank $3.$ Moreover, $\E$ is generated by elements $e_1, e_2, e_3$ belonging to $\Ecenter.$ The space of two forms $\twoform$ is isomorphic to $\A \oplus \A \oplus \A.$ The differential calculus $ ( \E, d ) $ is tame.  	
				\eppsn
				The tameness of the differential calculus $(\E, d )$ is observed in the proof of Theorem 6.6 of \cite{article1}.
				
			Let us  fix a torsionless connection $\nabla_0$ on $(\E, d )$ which we will need later.	Since $ \E: = \oneform $ is a free module with generators $ e_1, e_2, e_3, $ any connection on $ \E $ is determined by its action on $ e_1, e_2, e_3. $ Our choice of the torsion-less connection $ \nabla_0 $ is given by the following:
\be \label{27thjan2} \nabla_0 ( e_j ) = 0 ~ {\rm for} ~ j = 1,2; \nabla_0 ( e_3 ) = - e_1 \tensora e_2. \ee

The proof of the following proposition is a verbatim adaptation of the proof of Proposition 31 of \cite{chak_sinha} with the only difference that we use right connections instead of left connections. 
	\bppsn \label{30thmarch20n2}
 $ \nabla_0 $ is a torsion-less connection on $ \E. $
\eppsn

Now we define a pseudo-Riemannian bilinear metric $g_0$ on $ ( \E, d ). $ We will compute the Christoffel symbols and scalar curvature of the Levi-Civita connection for the triplet $( \E, d, g_0 )$ in 
	 Theorem \ref{explicitform}. 
		
		\blmma \label{24thjan5}
		 Let $ g_0 $ be the pseudo-Riemannian bilinear metric of Lemma \ref{27thmarch20} so that $ g_0 ( e_i \tensora e_j ) = \delta_{ij} .$ If $\nabla_0$ is the torsionless connection of Proposition \ref{30thmarch20n2} and $\Pi_g$ the map as in \eqref{20thfeb2020}, then  we have the following:
		  $$ \Pi_{g_0} ( \nabla_0 ) (  e_i \tensora e_j ) = - \sum_m  e_m T^m_{ij}, $$
			$$ {\rm where}, ~ T^2_{13} = T^2_{31} = 1 ~ {\rm and} ~ T^m_{ij} = 0 ~ {\rm otherwise}. $$
		\elmma
 {\bf Proof:} Let us begin by remarking that it is easy to see ( from Proposition \ref{30thmarch20n} ) that $(\E, d )$ satisfies the hypotheses of Lemma \ref{27thmarch20} so that the pseudo-Riemannian metric $g_0$ makes sense. Secondly, since $ \Pi_{g_0} ( \nabla_0 ) ( e_i \tensora e_j ) \in \E $ and $\E$ is a free right $\A$-module with basis $e_1, e_2, e_3,$ the elements $T^m_{ij}$ exist uniquely. 

 Since $ e_1, e_2, e_3 \in \Ecenter,$ we get
	    $$ \Pi_{g_0} ( \nabla_0 ) ( e_i \tensora e_j ) = \Pi^0_{g_0} ( \nabla_0 ) ( e_i \tensora e_j ) = ( g_0 \tensora {\rm id} ) \sigma_{23} ( \nabla_0 ( e_i ) \tensora e_j + \nabla_0 ( e_j ) \tensora e_i  ) = \Pi_{g_0} ( \nabla_0 ) ( e_j \tensora e_i ) $$
			by Proposition \ref{30thmarch203}.
			\be \label{24thjan4} {\rm Clearly},~ T^m_{ij} = T^m_{ji}. \ee
			 From \eqref{27thjan2}, it is immediate that for all $ i,j \in \{ 1, 2 \}, $
			$$ \Pi_{g_0} ( \nabla_0 ) ( e_i \tensora e_j ) = 0.  $$
			Moreover, 
			 \begin{eqnarray*}
			  \Pi_{g_0} ( \nabla_0 ) ( e_1 \tensora e_3 ) &=&  ( g_0 \tensora {\rm id} ) \sigma_{23} ( \nabla_0 ( e_1 ) \tensora e_3 + \nabla_0 ( e_3 ) \tensora e_1  )\\
				&=& -  ( g_0 \tensora {\rm id} ) \sigma_{23} (  ( e_1 \tensora e_2 ) \tensora e_1 )\\
				&=& -   ( g_0 \tensora {\rm id}  ) ( e_1 \tensora e_1 \tensora e_2 ) ~ {\rm (} ~ {\rm by} ~ \eqref{10thjuly20182} ~ {\rm )}\\
				&=& -  e_2.
			 \end{eqnarray*}
			Thus, 
			$$ T^2_{13} =  1, ~ T^1_{13} = T^3_{13} = 0. $$
			Next, 
			\begin{eqnarray*}
			  \Pi_{g_0} ( \nabla_0 ) ( e_3 \tensora e_2 ) &=&  ( g_0 \tensora {\rm id} ) \sigma_{23} ( \nabla_0 ( e_2 ) \tensora e_3 + \nabla_0 ( e_3 ) \tensora e_2  )\\
				&=& -  ( g_0 \tensora {\rm id} ) \sigma_{23} (  ( e_1 \tensora e_2 ) \tensora e_2 )\\
				&=& 0
			\end{eqnarray*}
			by \eqref{10thjuly20182} and as $ g ( e_i \tensora e_j ) = \delta_{ij}.1. $
			
			Thus, for all $ m = 1,2,3, ~ T^m_{32} = 0. $
			The rest of the $ T^m_{ij} $ can be computed by using \eqref{24thjan4}.
			\qed
			
			\brmrk
			A candidate $g^\prime$ for a canonical pseudo-Riemannian bilinear metric on a tame spectral triple was constructed in Subsection 3.1 of \cite{article1}. Proposition 6.4 of \cite{article1} verifies that $g^\prime$ satisfies all the required conditions to be a bilinear pseudo-Riemannian metric and moreover, the proof of this result shows that in fact $g^\prime$ satisfies \eqref{27thmarch202}. By the uniqueness of Lemma \ref{27thmarch20}, it follows that $g_0 = g^\prime.$     
			\ermrk
			
			Now, we are ready to compute the explicit form  of the Levi-Civita connection for the metric $ g_0 $ on the module $ \E. $ 
			
			\bthm \label{explicitform}
			Consider the tame diferential calculus $(\E, d )$ on $\A$ as above, $g_0$ be the pseudo-Riemannian bilinear metric of Lemma \ref{24thjan5} and $\nabla_0$ be the torsionless connection of Proposition \ref{30thmarch20n2}. 
			 Then the unique Levi-Civita connection $ \nabla $ for the triplet $ ( \E, d, g_0 ) $ is given by 
			$$ \nabla = \nabla_0 + L, $$ 
			where   $ L: \E \rightarrow \E \tensora \E $ is  defined by 
			\begin{equation} \label{30thmarch20n*} L ( e_j  ) = \sum_{i,m} e_i \tensora e_m L^j_{im}, \end{equation} 
			\be \label{28thfeb2017**} L^j_{im} =  \frac{1}{2} ( T^m_{ij} + T^i_{jm} - T^j_{mi} ), \ee
			where the elements $ \{ T^m_{ij} :  i,j,m = 1,2,3 \}  $ are as in Lemma \ref{24thjan5}.
			
			More precisely, the non zero $ L^j_{im} $ are as follows:
			$$ L^1_{23} = L^1_{32} = 0.5, ~ L^2_{13} = L^2_{31} = - 0.5, L^3_{12} = L^3_{21} = 0.5, $$
			where we have denoted $ \lambda 1_\A $ simply by $ \lambda. $
			If $ \nabla $ is given by 
			\be \label{28thfeb2017*} \nabla ( e_i ) = \sum_{j,k} e_j \tensora e_k \Gamma^i_{jk}, \ee
			then the non zero $ \Gamma^i_{jk} $ are as follows:
			$$ \Gamma^1_{12} = 1, \Gamma^1_{23} = \Gamma^1_{32} = 0.5, \Gamma^2_{12} = 1,  \Gamma^2_{13} = - 0.5, \Gamma^2_{31} = - 0.5, \Gamma^3_{12} = 1.5, \Gamma^3_{21} = 0.5. $$
			 \ethm
			{\bf Proof:} Theorem 6.6 of \cite{article1} proves that $(\E, d )$ is tame. By Remark \ref{22ndfeb2020}, we know that the map $ \Phi_{g_0} $ defined in Theorem \ref{19thfeb2020} is a right $\A$-linear isomorphism. Therefore, by Theorem \ref{19thfeb2020}, there exists a unique Levi-Civita connection $\nabla$ for the triplet $(\E, d, g_0 )$ and is given by
			$$ \nabla = \nabla_0 + \Phi^{-1}_{g_0} ( dg_0 - \Pi_{g_0} ( \nabla_0 ) ). $$
		Let us define $L = \nabla - \nabla_0.$ Then 	for all $i,j$ and for all $ X $ in $\E \tensorsym \E,$
			 \be \label{28thfeb20171} \Phi_{g_0} ( \nabla - \nabla_0 ) ( X ) = d g_0 ( X ) - \Pi_{g_0} ( \nabla_0 ) ( X ). \ee
			Morover, as $ \nabla $ and $ \nabla_0 $ are both torsion-less connections, we get
				\begin{equation} \label{28thfeb20172} \wedge L ( e_i ) = \wedge \nabla ( e_i ) - \wedge \nabla_0 ( e_i ) = - d ( e_i ) + d ( e_i ) = 0. \end{equation}
				
	Let $Q$ be the isomorphism from $ {\rm Ran} ( 1 - \Psym ) $ to $\twoform$ as in Definition \ref{defnQ}. Then
				\begin{eqnarray*}
				0 &=& Q^{-1} ( \wedge \circ L ( e_j )  ) ~ {\rm (} ~ {\rm by} ~ \eqref{28thfeb20172} ~ {\rm )}\\
				  &=& ( 1 - \Psym )  L ( e_j ) ~ {\rm (} ~ {\rm by} ~ \eqref{22ndmarch205} ~ {\rm )}\\
					&=& \frac{1 - \sigma}{2} ( \sum_{i,m} e_i \tensora e_m L^j_{im}  ) ~ {\rm (} ~ {\rm by} ~ \eqref{30thmarch20n*} ~ {\rm )}\\
					&=& \frac{1}{2} \sum_{i,m} e_i \tensora e_m ( L^j_{im} - L^j_{mi}  ) ~ {\rm (} ~ {\rm by} ~ \eqref{10thjuly20182} ~ {\rm )}.
					\end{eqnarray*}
	Since $ \{ e_i \tensora e_j ; i, j \} $ is a basis of $\E \tensora \E,$ we obtain			 
					  \be \label{28thfeb20173} L^j_{im} = L^j_{mi} ~ \forall ~ i,j,m. \ee
			Now we derive a relation between $ T^m_{ij} $ and $L^i_{jm}$ by the following computation:
			{\allowdisplaybreaks
				\begin{eqnarray*}
				&& \sum_m e_m T^m_{ij}\\
				&=& \frac{1}{2} \sum_m e_m T^m_{ij} + \frac{1}{2} \sum_m e_m T^m_{ji} ~ {\rm (} ~ {\rm by} ~ \eqref{24thjan4} ~ {\rm )}\\
				&=& - \frac{1}{2} \Pi_{g_0} ( \nabla_0 ) ( e_i \tensora e_j ) - \frac{1}{2} \Pi_{g_0} ( \nabla_0 ) ( e_j \tensora e_i ) ~ {\rm (} ~ {\rm Lemma} ~ \ref{24thjan5} ~ {\rm )}\\
				&=& - \Pi_{g_0} ( \nabla_0 ) ( \frac{e_i \tensora e_j + e_j \tensora e_i}{2}  )\\
				&=& d ( g_0 ( \frac{e_i \tensora e_j + e_j \tensora e_i}{2} ) ) - \Pi_{g_0} ( \nabla_0 ) ( \frac{e_i \tensora e_j + e_j \tensora e_i}{2}   ) ~ {\rm (} {\rm since} ~ g_0 ( e_i \tensora e_j ) = \delta_{ij}.1_{\A} ~ {\rm )}\\
				&=& \Phi_{g_0} ( L ) ( \frac{e_i \tensora e_j + e_j \tensora e_i }{2}  )\\
				&& {\rm (} ~ {\rm by} ~ \eqref{28thfeb20171} ~ {\rm and} ~ {\rm since} ~ \frac{e_i \tensora e_j + e_j \tensora e_i}{2} = \Psym ( e_i \tensora e_j ) \in {\rm Ran} ( \Psym ) = \E \tensorsym \E ~ {\rm )}\\
				&=& ( g_0 \tensora {\rm id}  ) \sigma_{23} ( L \tensora {\rm id} ) ( 1 + \sigma ) ( \frac{e_i \tensora e_j + e_j \tensora e_i}{2} ) ~ {\rm (} ~  {\rm by} ~ {\rm the} ~ {\rm definition} ~ {\rm of} ~ \Phi_{g_0} ~ {\rm in} ~ {\rm Theorem} ~ \ref{19thfeb2020} ~ {\rm )}\\
				&=& ( g_0 \tensora {\rm id}  ) \sigma_{23} ( L \tensora {\rm id} ) ( 1 + \sigma ) ( \frac{1 + \sigma}{2} ) ( e_i \tensora e_j ) ~ {\rm (} ~  {\rm by}  ~ \eqref{10thjuly20182} ~ {\rm )}\\
				&=& 2 ( g_0 \tensora {\rm id}  ) \sigma_{23} ( L \tensora {\rm id} ) ( \frac{1 + \sigma}{2} )^2 ( e_i \tensora e_j )\\
				&=& 2 ( g_0 \tensora {\rm id}  ) \sigma_{23} ( L \tensora {\rm id} ) \Psym ( e_i \tensora e_j ) ~ {\rm (} ~ {\rm as} ~ \Psym = \frac{1 + \sigma}{2} ~ {\rm is} ~ {\rm an} ~ {\rm idempotent} ~ {\rm )}\\
				&=& ( g_0 \tensora {\rm id}  ) \sigma_{23} ( L \tensora {\rm id} ) ( e_i \tensora e_j + e_j \tensora e_i )\\
				&=& ( g_0 \tensora {\rm id}  ) \sigma_{23} ( \sum_{k,x} e_k \tensora e_x L^i_{kx} \tensora e_j + \sum_{l,y} e_l \tensora e_y L^j_{ly} \tensora e_i   ) ~ {\rm (} ~ {\rm by} ~ \eqref{30thmarch20n*} ~ {\rm )}\\
				&=& ( g_0 \tensora {\rm id} ) ( \sum_{k,x} e_k \tensora e_j \tensora e_x L^i_{kx} + \sum_{l,y} e_l \tensora e_i \tensora e_y L^j_{ly} ) ~ {\rm (} ~ {\rm by} ~ \eqref{10thjuly20182} ~ {\rm )}\\
				&=& \sum_k \delta_{jk} \sum_{x} e_x L^i_{kx} + \sum_l \delta_{li} \sum_y e_y L^j_{ly}\\
				 &=& \sum_y e_y ( L^i_{jy} + L^j_{iy} ).
				\end{eqnarray*}
			}
				This implies that
				\be \label{24thjan6} L^i_{j,m} + L^j_{i,m} = T^m_{ij}. \ee
				Interchanging $ ( i,j,m  ) $ with $ ( j,m,i  ) $ and $ ( m,i,j  ) $, we have respectively:
				\be \label{24thjan7} L^j_{m,i} + L^m_{j,i} = T^i_{j,m}, \ee
				 \be \label{24thjan8} L^m_{i,j} + L^i_{m,j} = T^j_{m,i}. \ee
				Now, by \eqref{24thjan6} + \eqref{24thjan7} - \eqref{24thjan8} and \eqref{28thfeb20173}, we have
				$$L^j_{i,m} = \frac{1}{2} ( T^m_{ij} + T^i_{jm} - T^j_{mi}  ), $$
				which proves \eqref{28thfeb2017**}.
				The numerical expressions for $ L^{i}_{j,m}  $ follow from the values of $T^m_{ij}$ in Lemma \ref{24thjan5}. Finally, since $ \nabla = \nabla_0 + L, $ the Christoffel symbols $ \Gamma^i_{jk} $ as in \eqref{28thfeb2017*} can be computed by using \eqref{27thjan2} and the values of $L^i_{j,m}.$
			 \qed

						\bthm \label{theoremfinal}
						Let $ \nabla $ denote the Levi-Civita connection for the metric $ g_0 $ on the module $ \E $ of one forms over the quantum Heisenberg manifold $ \cla. $ The the  Ricci and scalar curvature of $ \nabla $ are as follows:
						$$ {\rm Ric} ( e_1, e_1 ) = - 1, ~ {\rm Ric} ( e_2, e_2 ) = 1, {\rm Ric} ( e_1, e_3 ) =  {\rm Ric} ( e_3, e_3 ) = - 0.5, {\rm Ric} ( e_2, e_3 ) = - 0.5,  $$
						$$ {\rm Ric} ( e_1, e_2 )  =  {\rm Ric} ( e_2, e_1 )   = {\rm Ric} ( e_3, e_1 ) = {\rm Ric} ( e_3, e_2 ) = 0. $$
						$${\rm Scal} = - 0.5.$$
						\ethm
						{\bf Proof:} The proof follows by a direct computation using Theorem \ref{explicitform} and  the formulas  of $ {\rm Ric} ( e_j, e_l ) $ and $ {\rm Scal} $  in Proposition \ref{20thjanprop0}. \qed
						
						\brmrk
						 We note that the quantum Heisenberg manifold has a constant negative scalar curvature and moreover, the curvature is independent of the choice of the parameter $ \alpha $ used to define the Dirac operator. 
						\ermrk

{\bf Acknowledgement}  The authors would  like to thank Aritra Bhowmick for helping them with the computer aided computations of Theorem \ref{explicitform} and Theorem \ref{theoremfinal}.  D.G will like to thank D.S.T, Government of India for J.C. Bose National Fellowship. S.J will like to thank D.S.T, Government of India for the Inspire Fellowship.


\begin{thebibliography}{6666}

\bibitem{pseudo} J. Arnlind and M. Wilson: Riemannian curvature of the noncommutative 3-sphere, J. Noncomm. Geom. {\bf 11} (2017)  507--536.

\bibitem{cylinder} J. Arnlind and G. Landi:  Projections, modules and connections for the noncommutative cylinder, Adv. Theor. Math. Phys. 24 (2020) 527--562.

\bibitem{tiger} J. Arnlind and A.T. Norkvist: Noncommutative minimal embeddings and morphisms of pseudo-Riemannian calculi, J. Geom. Phys. {\bf 159} (2020).

\bibitem{landiqhom} J. Arnlind, K. Ilwale and G. Landi: On q-deformed Levi-Civita connections, arXiv: 2005.02603.

\bibitem{aschieri} P. Aschieri: Cartan structure equations and Levi-Civita connection in braided geometry, arXiv:2006.02761.


\bibitem{majid_2} E.J. Beggs and S. Majid: $\ast$-Compatible connections in noncommutative Riemannian geometry, J. Geom. Phys. {\bf 61} (2011) 95--124.


\bibitem{beggsbrezinski} E.J. Beggs and T. Brzezinski: Noncommutative differential operators, Sobolev spaces and the centre of a category, Journal of Pure and Applied Algebra {\bf 218} ( 2014 ) 1--17.

\bibitem{chern} E.J. Beggs and S. Majid: Spectral triples from bimodule connections and Chern connections, Journal of Noncommutative Geometry, {\bf 11} (2) ( 2017 ), 669--701.

\bibitem{beggsmajidbook} E.J. Beggs and  S. Majid: Quantum Riemannian geometry, Grundlehren der mathematischen Wissenschaften, Springer Verlag, 2019. 


\bibitem{article1} J. Bhowmick, D. Goswami and S. Mukhopadhyay: Levi-Civita connections for a class of spectral triples, Letters in Mathematical Physics, {\bf 110} ( 2020 ), 835--884. 
			
			\bibitem{article3} J. Bhowmick, D. Goswami and G. Landi: On the Koszul formula in noncommutative geometry, {Reviews in Mathematical Physics}, {\bf 32}, No 10, 2050032, 2020.

			
\bibitem{article4} J. Bhowmick, D. Goswami and G. Landi: Levi-Civita connections and vector fields for noncommutative differential calculi, {Internat. J. Math.}, {\bf 31}, No 8, 2020.
			
		
			
			\bibitem{article2} J. Bhowmick, D. Goswami and S. Joardar: A new look at Levi-Civita connection in noncommutative geometry, arXiv: 1606.08142v5.

			
					
		\bibitem{sitarz}	A. Bochniak, A. Sitarz and P. Zalecki: Riemannian geometry of a discretized circle and torus, arXiv:2007.01241.		
					
	\bibitem{chak_sinha} P. S. Chakraborty and K.B. Sinha: Geometry on the Quantum Heisenberg Manifold, Journal of Functional Analysis, {\bf 203}, 425--452, 2003.


 \bibitem{connes} A. Connes: Noncommutative geometry. Academic Press, San Diego, CA, 1994.
			
			
			\bibitem{Connes-dubois} A. Connes and M. Dubois-Violette : Noncommutative finite-dimensional manifolds. I.,
Spherical manifolds and related examples. Comm.\ Math.\ Phys. \textbf{230} (2002), no. 3, 539--579.


\bibitem{connes_landi} A. Connes and G. Landi: Noncommutative Manifolds the Instanton Algebra and Isospectral Deformations, Commun. Math. Phys.221:141--159,2001.


\bibitem{Connes_moscovici} A. Connes and H. Moscovici: Modular curvature for noncommutative two-tori, 
J. Amer. Math. Soc. {\bf 27 }(2014) 639--684.

\bibitem{scalar_3} A. Connes and P. Tretkoff: The Gauss-Bonnet theorem for the noncommutative two
torus, in Noncommutative Geometry, Arithmetic, and Related Topics, Johns Hopkins Univ. Press, Baltimore, MD, 2011, 141--158.

\bibitem{dubois} M. Dubois-Violette and P.W. Michor: Derivation et calcul differentiel non commutatif II, C. R. Acad. Sci. Paris Ser. I Math, {\bf 319} (1994) 927--931. 

\bibitem{dubois2} M. Dubois-Violette and P.W. Michor: Connections on central bimodules, J. Geom. Phys. {\bf 20} (1996) 218--232.


\bibitem{khalkhali} F. Fathizadeh, M. Khalkhali: Curvature in Noncommutative Geometry, arXiv:1901.07438.

 
\bibitem{frolich} J. Frohlich, O. Grandjean and A. Recknagel: Supersymmetric Quantum
Theory and Non-Commutative Geometry, Commun. Math. Phys {\bf 203} (1999) 119--184.

	
	\bibitem{heckenberger_etal}   I. Heckenberger and K. Schmuedgen : Levi-Civita Connections on the Quantum Groups 
	${\rm SL}_q(N)$ ${\rm O}_q(N)$ and ${\rm Sp}_q(N)$,
Comm. Math. Phys. {\bf 185} (1997) 177--196.

\bibitem{soumalya} S. Joardar: Scalar Curvature of a Levi-Civita Connection on Cuntz algebra with three generators, Lett. Math. Phys., {\bf 109}, 2665--2679, 2019. 


\bibitem{majid_1} S. Majid: Noncommutative Riemannian and spin geometry of the standard q-sphere, Commun.
Math. Phys. {\bf 256 }(2005) 255--285.

\bibitem{matassa} M. Matassa: Fubini-Study metrics and Levi-Civita connections on quantum projective spaces, arXiv: 2010. 03291. 

\bibitem{sheu} M.A. Peterka and A.J.L. Sheu:  On Noncommutative Levi-Civita Connections, International Journal of Geometric Methods in Modern Physics {\bf 14}, No. 5 (2017) 1750071.

\bibitem{rieffel_heisenberg}  M.A. Rieffel: Deformation quantization of Heisenberg manifolds, \emph{Comm.~ Math.~ Phys.} \textbf{122} (1989),  531--562.

\bibitem{rieffel}  M. A. Rieffel : Deformation Quantization for actions of $ \IR^{d},$ Memoirs of the American Mathematical Society, November 1993. Volume {\bf 106}, Number 506.

\bibitem{Rosenberg} J. Rosenberg: Levi-Civita's Theorem for Noncommutative Tori, SIGMA {\bf 9 }(2013),
071.

\bibitem{weber} T. Weber: Braided Cartan Calculi and Submanifold Algebras, arXiv: 1907.13609.


\bibitem{Skd}  M. Skeide: Hilbert modules in quantum electrodynamics and quantum probability.  Comm. Math. Phys. \textbf{192} (1998), 569--604.


   
\end{thebibliography}
\end{document}